\documentclass[10pt,peerreview]{IEEEtran}
\usepackage{amsmath}
\usepackage{amssymb}
\usepackage{bm}
\usepackage[dvipsnames]{xcolor}
\usepackage{hyperref}

\usepackage{tikz}
\usetikzlibrary{calc}
\usetikzlibrary{shapes.geometric}
\usetikzlibrary{decorations.markings}
\usetikzlibrary{decorations.pathmorphing, patterns,shapes,arrows}
\usetikzlibrary{automata}
\usetikzlibrary{shadows}
\usetikzlibrary{shadings}
\usetikzlibrary{decorations.pathreplacing}

\newtheorem{lemma}{Lemma}[section]
\newtheorem{prop}{Proposition}[section]
\newtheorem{dfn}{Definition}[section]

\newtheorem{exm}{Example}[section]
\newtheorem{rem}{Remark}[section]
\newtheorem{hyp}{Hypothesis}[section]
\numberwithin{equation}{section}
\def\theequation{\arabic{section}.\arabic{equation}}

\newcommand{\RgO}{\ensuremath{\mathbb{R}_{>0}}}
\newcommand{\RgeO}{\ensuremath{\mathbb{R}_{ \geq 0}}}
\newcommand{\Rat}[1]{\ensuremath{\mathbb{R}^{#1}}}
\newcommand{\U}{\mathcal{U}}

\begin{document}
\title{Observer design for triangular systems under weak observability assumptions}

	\author{D.~Theodosis,  D.~Boskos and J.~Tsinias %
		\thanks{D. Theodosis and J. Tsinias are with the Department of Mathematics, National Technical University of Athens, Zografou Campus 15780, Athens, Greece, email: dtheodp@central.ntua.gr (corresponding author), jtsin@central.ntua.gr.}%
		\thanks{D. Boskos is with the ACCESS Linnaeus Centre, School of Electrical Engineering, KTH Royal Institute of Technology, SE-100 44, Stockholm, Sweden and with the KTH Centre for Autonomous Systems, email: boskos@kth.se}} %

\maketitle

\begin{abstract}
This paper presents results on the solvability of the observer design problem for general nonlinear triangular systems with inputs, under weak observability assumptions. {The local state estimation is exhibited by means of a delayed time-varying Luenberger-type system. In order to achieve the global estimation, a switching sequence of observers is designed.}
\end{abstract}

\begin{IEEEkeywords}observer design, nonlinear triangular systems, switching dynamics.\end{IEEEkeywords}

\section{Introduction}

Observer design for nonlinear systems constitutes a central problem in control theory with several contributions during the last decades; see for instance  \cite{AAtLLf99}-\cite{TjTd17}. Several approaches have been leveraged for the solvability of this problem, including high-gain, Lyapunov-based, switching and various other techniques (see e.g.,  \cite{AaRa15}, \cite{AvPlAa09}, \cite{BgTa07}, \cite{GjHhOs92}, \cite{GjKi01}, \cite{HhTbAf02}, \cite{KpKf04};  \cite{AjKh09}, \cite{AaPl06}, \cite{DhQcYsLs13}, \cite{TJ08}; \cite{AjKh09}, \cite{BdTj13a}, \cite{HPR14}, \cite{Liu97}; \cite{AAtLLf99}, \cite{Am99}, \cite{AvPl06},  \cite{BjSjh08}, \cite{Jp03}-\cite{KajXm02}, \cite{Pl03}).

In the present work we derive sufficient conditions for the solvability of the observer design problem (ODP) for a class of time-varying nonlinear triangular control systems of the form\begin{subequations} \label{system:triangular}
	\begin{align}
	\dot{x}_{i}=&f_{i}(t,x_{1},\ldots,x_{i},u)+a_{i}(t,x_{1},u)x_{i+1},\, \nonumber \\
		&i=1,2,\ldots,n-1; \label{system:triangular:st} \\
		\dot{x}_{n}=&f_{n} (t,x_{1} ,\ldots,x_{n},u),  \nonumber \\
	y=&x_{1}, (x_{1} ,\ldots,x_{n})\in\Rat{n} ,u\in\mathbb{R}^{p}\label{system:triangular:out}
	\end{align}
\end{subequations}

\noindent where $u(\cdot)$, $y(\cdot)$ are the input and output of the system, respectively. There is an extensive literature concerning the design of observers for triangular systems. For instance, in \cite{AvPlAa09}, the time-invariant case of \eqref{system:triangular} with measurement noise is considered under the assumption that  $a_{i}(\cdot)\ge \rho$, $i=1,\ldots,n-1$ for some constant  $\rho>0$. The authors presented a new class of high-gain observers with updating gain which extends the classical result in \cite{GjKi01}. Under the same assumption, the high-gain ODP is explored in \cite{KpKf04}, where the authors propose a dynamic high-gain scaling technique. In \cite{DhQcYsLs13}, sufficient conditions for the existence of finite-time convergent observers are provided for time-varying triangular systems \eqref{system:triangular}, under a similar assumption for the functions $a_i$. There are also several works where the ODP is explored for a class of systems \eqref{system:triangular} with $a_{i}\equiv 1$, $i=1,\ldots,n-1$; see for instance \cite{AAtLLf99}, \cite{AaRa15}, \cite{AaPl06}, \cite{BjSjh08}, \cite{GjHhOs92}, \cite{HhTbAf02}, \cite{Jp03} and \cite{Pl03}. Switching techniques have been used in \cite{AjKh09}, \cite{BdTj13a}, \cite{HPR14} and \cite{Liu97} for the solvability of the state estimation for certain class of nonlinear systems, under appropriate hypotheses.

The main contribution of this paper is to extend and generalize previous authors work \cite{BdTj13a}, concerning the solvability of the ODP for triangular control systems of the form \eqref{system:triangular}, by means of a \textit{delayed switching observer}.
The paper is organized as follows. Section \ref{section:preliminaries} contains the notations and various concepts, including the concept of the switching observer for general time-varying systems:
\begin{subequations} \label{system:nonlinear}
\begin{align}
\dot{x}&=f(t,x,u),(t,x,u)\in {\mathbb R}_{\ge 0} \times {\mathbb R}^{n}\times\Rat{p} \label{system:nonlinear:st} \\
y&=h(t,x,u),y\in {\mathbb R}^{k} \label{system:nonlinear:out}
\end{align}
\end{subequations}

\noindent where $y(\cdot )$ is the output and $u(\cdot)$ is the input of the system. We then provide the precise statement of our main result ({Proposition} \ref{Theorem}) concerning the solvability of the ODP for \eqref{system:triangular}. Section \ref{section:Observers:general_sys} contains some preliminary results concerning solvability of the delayed ODP for the general case \eqref{system:nonlinear} with linear output (Propositions \ref{Proposition:state_deter} and \ref{Proposition:switching}). Finally, in Section \ref{section:main_result}, we use the results of Section \ref{section:Observers:general_sys}, in order to prove our main result.
\section{Notations, Definitions and Main Result }\label{section:preliminaries}

\subsection{Notations} Throughout this paper we adopt the following notation.
For a given vector $x\in{\mathbb R}^{n}$, $x'$ denotes its transpose and $|x|$ its Euclidean norm.
We use the notation $|A|:=\max \{|Ax|:x\in {\mathbb R}^{n};|x|=1\}$ for the induced norm of a matrix $A\in {\mathbb R}^{m\times n}$ and $|A|_{F}$ for its Frobenius norm, namely $|A|_{F}=\left(\sum_{i=1}^{m}\sum_{i=1}^{n}a_{i,j}^{2}\right)^{\frac{1}{2}}$. We denote by ${\rm diag}\{a_1,\ldots,a_n\}$ the diagonal matrix with entries $a_1,\ldots,a_n\in\Rat{}$ and by $I_{n\times n}$ the $n\times n$ identity matrix.
By $N$ we denote the class of all increasing $C^{0} $ functions $\phi:{\mathbb R}_{\ge 0}\to{\mathbb R}_{\ge 0} $.
For given $R>0$, we denote by $B_{R}$ the closed ball of radius $R>0$, centered at $0\in{\mathbb R}^{n} $.
Consider a pair of metric spaces $X_{1}$, $X_{2}$ and a set-valued map $X_{1}\ni x\to Q(x)\subset X_{2}$. We say that $Q(\cdot)$ satisfies the \textit{Compactness Property }(\textbf{CP)}, if for every sequence $(x_{\nu})_{\nu \in {\mathbb N}}\subset X_{1}$ and $(q_{\nu})_{\nu \in {\mathbb N}}\subset X_{2}$ with $x_{\nu }\to x\in X_{1}$ and $q_{\nu} \in Q(x_{\nu})$, there exist a subsequence $(x_{\nu_{k}} )_{k\in{\mathbb N}}$ and $q\in Q(x)$ such that $q_{\nu_{k}}\to q$.
{Given $t_0\ge 0$, $\tau>0$, a nonempty set $S$ and a function $g:[t_0,\infty)\to S$, we define its $\tau$-time shift $g_{\tau}:[t_0+\tau,\infty)\to S$ as $g_{\tau}(t):=g(t-\tau)$, $t\in [t_0+\tau,\infty)$.}

\subsection{Definitions and Main Result}

We assume that the right hand side of \eqref{system:nonlinear:st} is locally Lipschitz with respect to $x$, i.e., for each compact $I\subset \RgeO$, $K\subset \Rat{n}$ and $U\subset \Rat{p}$ there exists a constant $L>0$ such that
$$
|f(t,x,u)-f(t,z,u)|\le L|x-z|,\forall t\in I, x,z\in K, u\in U.
$$

 We next provide the definition of $(M,\U)$-forward completeness for system \eqref{system:nonlinear}, which constitutes a generalization of the classical forward completeness property.

\begin{dfn}\label{dfn:forw:comp}
Consider a nonempty subset $M$ of $\Rat{n}$. For each $(t_0,x_0)\in\RgeO\times M$, let $\U(t_0,x_0)$ be a nonempty set of (measurable and locally essentially bounded) inputs $u:[t_0,\infty)\to \Rat{p}$ and define
\begin{equation} \label{input:set}
\U(t_0):=\cup_{x_0\in M}\U(t_0,x_0),t_0\ge 0;\; \U:=\cup_{t_0\in\RgeO}\U(t_0).
\end{equation}

\noindent We say that system \eqref{system:nonlinear:st} is $(M,\U)$-forward complete, if there exists a function $\beta\in NN$  such that the solution $x(\cdot):=x(\cdot,t_{0},x_{0};u)$ of \eqref{system:nonlinear:st} corresponding to input ${u(\cdot)}$ and initiated from $x_0$ at time $t=t_{0} $ {is defined for all $t\ge t_0$ and} satisfies
\begin{equation} \label{state:bound}
|x(t)|\le\beta (t,|x_{0}|), \forall t\ge t_{0}\ge 0, x_{0} \in M,u\in\U(t_0,x_0).
\end{equation}
\end{dfn}


It turns out, that under $(M,\U)$-forward completeness of \eqref{system:nonlinear:st}, for each $t_{0}\ge 0$, $x_{0}\in M$ and $u\in\U(t_0,x_0)$, the corresponding output $y(t)=h(t,x(t,t_{0},x_{0};u),u(t))$ of \eqref{system:nonlinear} is defined for all $t\ge t_{0}$. For each $t_{0} \ge 0$ and $x_{0} \in M$ we consider the set $O(t_{0},x_0)$, containing the pairs of all possible inputs in $\U(t_0,x_0)$ and their corresponding output paths of system \eqref{system:nonlinear} initiated at $(t_0,x_0)$, namely:
\begin{align*}
O(t_{0},x_0):=\{&(u,y):[t_{0},\infty)\to {\mathbb R}^{p}\times\mathbb{R}^{k}: u\in \U(t_0,x_0), y(t)=h(t,x(t,t_{0},x_{0}{;u),u(t)}),\forall t\ge t_0\}.
\end{align*}

\noindent Define
\begin{equation} \label{output:functions}
O(t_{0},M):=\cup_{x_0\in M} O(t_{0},x_0).
\end{equation}


\begin{dfn}\label{dfn:causal}
Let $k,\ell,m,n\in\mathbb{N}$, $\emptyset\ne M\subset\Rat{n}$, $\emptyset\ne S\subset\Rat{\ell}$ and for each $t_{0}\ge 0$, let $\Omega(t_{0},M)$ be a nonempty set of functions $(u,y):[t_{0},\infty)\to\Rat{p}\times \Rat{k}$, $y:=y_{t_0,x_0}$, $u:={u_{t_0,x_0}}$ parameterized by $t_0\ge 0$ and $x_{0}\in M$. Given $I\subset [t_0,\infty)$, we say that the map
\begin{equation*}
I\times\Omega(t_{0},M)\ni(t,(u,y))\to a_{u,y}(t)\in S
\end{equation*}

\noindent is causal with respect to $\Omega(t_{0},M)$, if for each $t\in I$, the value $a(t):=a_{u,y}(t)$ depends only on $(u,y)|_{[t_{0},t)}$ (the restriction of $(u,y)(\cdot)$ on $[t_{0},t)$). {Let $\beta>\alpha\ge t_0$ and $a_{u,y}(\cdot)$ as defined above. We say that $a_{u,y}(\cdot)$ is strongly causal on $I\cap [\alpha,\beta]$ with respect to $\Omega(t_0,M)$, if for each $t\in I\cap [\alpha,\beta]$, the value $a(t)(=a_{u,y}(t))$ depends only on $(u,y)|_{[t_0,\alpha)}$, namely, on the values of $(u,y)$ on the interval $[t_0,\alpha)$ and is generally independent of the values of $(u,y)$ on the interval $I\cap [\alpha,\beta]$.}

\end{dfn}


\indent An observer for the general deterministic system \eqref{system:nonlinear} is a system,  driven by both the input and output of \eqref{system:nonlinear}, which achieves the online estimation of the state of \eqref{system:nonlinear}. In this work, we investigate the ODP for system \eqref{system:triangular}, under the hypothesis that each $a_{i}(\cdot,y(\cdot),u(\cdot))$  may vanish on certain subintervals of $\RgeO$. In particular, at the initialization of the system, we assume knowledge of a partition of $\RgeO$ into a countable sequence of intervals, each of which containing an instant where all the $a_i$'s will be nonzero. However, there is \textit{no a priori knowledge} of these time instants. Thus, in order to construct the desired observer, we require some future knowledge of the output of the system, resulting in a delayed state estimation. This delay constitutes a design parameter which can be tuned arbitrarily small.

In order to formalize the approach discussed above, we first introduce the concept of the \textit{delayed} observer as well as the concept of the \textit{delayed switching observer} fo the general case \eqref{system:nonlinear}. 

\begin{dfn} \label{Mfc:def}
Let $\emptyset\ne M\subset\Rat{n}$, $\U$ as in \eqref{input:set} and assume that system \eqref{system:nonlinear} is $(M,\U)$-forward complete. Given $\tau>0$, we say that the \textbf{$\tau$-Delayed Observer Design Problem ($\tau$-DODP) is solvable for \eqref{system:nonlinear} with respect to $(M,\U)$}, if for every $t_{0} \ge 0$, and $(u,y)\in O(t_{0},M)$ there exist a continuous map
\begin{equation*}
G:=G_{t_{0},\tau,{y,u}}{(t,z,w,u)}:[t_{0}+\tau,\infty)\times {\mathbb R}^{n}\times {\mathbb R}^{k}\times\mathbb{R}^{p}\to {\mathbb R}^{n},
\end{equation*}

\noindent causal with respect to $O(t_{0},M)$ and a nonempty set $\bar{M}\subset {\mathbb R}^{n}$ such that for every $z_{0} \in \bar{M}$ the corresponding trajectory $z(\cdot):=z(\cdot,t_{0}+\tau,z_{0};u,y)$; $z(t_{0}+\tau)=z_{0}$ of the observer
\begin{equation*}
\dot{z}(t)=G(t,z(t),y_{\tau}(t),u_{\tau}(t))
\end{equation*}

\noindent exists for all $t\ge t_{0}+\tau $ and the error $e(t):={x_{\tau}(t)-z(t)}$, between the trajectory $x(\cdot):=x(\cdot,t_{0},x_{0};u)$, $x_{0} \in M$ of \eqref{system:nonlinear:st} and the trajectory $z(\cdot ):=z(\cdot,t_{0}+\tau,z_{0};u,y)$ of the observer satisfies:
\begin{equation} \label{error:conv}
\mathop{\lim }\limits_{t\to \infty } e(t)=0.
\end{equation}

\noindent We say that the \textbf{Infinitesimally Delayed Observer Design Problem (IDODP) is solvable for \eqref{system:nonlinear} with respect to} $(M,\U)$, if the $\tau$-DODP is solvable for \eqref{system:nonlinear} for any arbitrarily small $\tau>0$.  
\end{dfn}

\begin{dfn}\label{dfn:Dsods}
\noindent Let $\emptyset\ne M\subset\Rat{n}$, $\U$ as in \eqref{input:set} and assume that system \eqref{system:nonlinear} is $(M,\U)$-forward complete. Given $\tau>0$, we say that the \textbf{$\tau$-Delayed Switching Observer Design Problem ($\tau$-DSODP) is solvable for \eqref{system:nonlinear} with respect to  $(M,\U)$}, if for every $t_{0} \ge 0$ and $(u,y)\in O(t_{0},M)$ there exist a strictly increasing sequence of times $(t_{m})_{m\in {\mathbb N}}$ with
\begin{equation*}
t_{1}=t_{0}+\tau\;{\rm and}\;\lim_{m\to\infty} t_{m}=\infty,
\end{equation*}

\noindent a sequence of continuous mappings
\begin{equation*}
G_{m}:=G_{m,t_{m-1},\tau,u,y}(t,z,w,u):[t_{m-1},t_{m+1}]\times {\mathbb R}^{n}\times {\mathbb R}^{k}\times\mathbb{R}^{p}\to {\mathbb R}^{n},m\in {\mathbb N},
\end{equation*}

\noindent causal with respect to $O(t_{0},M)$, and a nonempty set $\bar{M}\subset {\mathbb R}^{n}$ such that the solution $z_{m}(\cdot)$ of the system
\begin{equation} \label{observer:seq}
\dot{z}_{m}(t)=G_{m}(t,z_{m}(t),y_{\tau}(t),u_{\tau}(t)), t\in [t_{m-1},t_{m+1}]
\end{equation}

\noindent with initial $z(t_{m-1})\in \bar{M}$, is defined for every $t\in [t_{m-1},t_{m+1}]$ and in such a way that, if we consider the piecewise continuous map $Z:[t_{0}+\tau,\infty )\to {\mathbb R}^{n} $ defined as $Z(t):=z_{m}(t)$, $t\in [t_{m},t_{m+1})$, $m\in {\mathbb N}$, where for each $m\in {\mathbb N}$, $z_{m}(\cdot)$ denotes the solution of \eqref{observer:seq}, then the error $e(t):=x_{\tau}(t)-Z(t)$ between the trajectory $x(\cdot):=x(\cdot,t_{0},x_{0};u)$, of \eqref{system:nonlinear:st} and $Z(\cdot)$ satisfies \eqref{error:conv}. We say that the \textbf{Infinitesimally Delayed Switching Observer Design Problem (IDSODP) is solvable for \eqref{system:nonlinear} with respect to} $(M,\U)$, if the $\tau$-DSODP is solvable for \eqref{system:nonlinear} for any $\tau>0$. 
\end{dfn}


{We provide now the precise statement of the main result of present work for the solvability of the IDOSDP (IDODP) for triangular systems \eqref{system:triangular}. We assume that for each $i=1,\ldots,n$ the map $f_i:\RgeO\times\Rat{i}\times\Rat{p}\to\Rat{}$ is $C^0$, for each fixed $t\ge 0$ and $u\in \Rat{p}$, $f_i(t,\cdot,u):\Rat{n}\to \Rat{}$ is $C^1$ and for every $i=1,\ldots,n-1$ the map $a_i:\RgeO\times\Rat{}\times\Rat{p}\to\Rat{}$ is $C^1$. Moreover we assume:}

\noindent \textbf{H1.} {There} exist a nonempty subset $M$ of $\Rat{n}$, a continuous function $\bar{u}:\RgeO\to\RgO$ and a nonempty set of \textbf{continuously differentiable} inputs $\U$ as in \eqref{input:set}, such that \eqref{system:triangular:st} is $(M,\U)$-forward complete; particularly, assume that there exists a function $\beta\in NN$ such that the solution $x(\cdot):=x(\cdot,t_0,x_0;u)$ of \eqref{system:triangular:st} satisfies \eqref{state:bound}. Additionally, we assume that for each $t_0\ge 0$ and $u\in\U(t_0)$ it holds
\begin{equation}\label{u:bound}
|u(t)|\leq \bar{u}(t),\forall t\geq t_0\ge 0.
\end{equation}

\noindent \textbf{H2.} For every $t_{0}\ge 0$ and $(u,y)\in O(t_0,M)$, there exists an \textit{a priori} known strictly increasing sequence of times $\{T_{\nu}\}_{\nu\in\mathbb{N}_{0}}$ with
\begin{equation}
T_{0}=t_{0};\lim_{\nu\to\infty}T_{\nu}=\infty
\end{equation}

\noindent in such a way that a sequence $\{\hat{t}_{\nu}\}_{\nu\in\mathbb{N}}$ can be found with $\hat{t}_{\nu}\in(T_{\nu-1},T_{\nu})$ for all $\nu\in\mathbb{N}$, such that
\begin{equation}\label{persistency:condition}
a_{i}(\hat{t}_{\nu},y(\hat{t}_{\nu}),u(\hat{t}_{\nu}))\ne 0,\forall i=1,\ldots,n-1,\nu\in\mathbb{N}.
\end{equation}

\noindent Our main result is the following proposition.

\begin{prop}\label{Theorem}
{For system \eqref{system:triangular}, assume that there exists a nonempty subset $M$ of $\Rat{n}$ and a set of inputs $\U$ as in \eqref{input:set} such that H1 and H2 are fulfilled. Then }

\noindent (i) the IDSODP is solvable for \eqref{system:triangular} with respect to $(M,\U)$.

\noindent (ii) if in addition we assume that it is a priori known, that the initial states of \eqref{system:triangular} belong to the (nonempty) intersection of $M$ with a given ball $B_{R}$ of radius $R>0$ centered at zero $0\in {\mathbb R}^{n}$, then the IDODP is solvable for \eqref{system:triangular} with respect to $(B_{R}\cap M,\U)$.
\end{prop}


The following elementary example illustrates the nature of Proposition \ref{Theorem}.
\begin{exm}\label{example:1}
Consider the system
\begin{equation}\label{exm:1:system}
\dot{x}_1=ux_2,\dot{x}_2=g(t,x_1,x_2,u)-x_2^q, y=x_1,\,\,(x_1,x_2)\in\Rat{2}
\end{equation}
		
\noindent {where $q$ is an odd integer and} $g\in C^1({\Rat{3};\Rat{})}$ satisfies 
\begin{equation}\label{function:g:property}
|g(t,x_1,x_2,u)|\le \alpha(t)|x_1|+\beta(t),\forall x_1,x_2\in\Rat{}, |u|\le \bar{u}(t), t\ge 0,
\end{equation}

\noindent for certain continuous function $\bar{u}:\RgeO\to\RgO$ and $\alpha$, $\beta\in N$.  We assume that the input set $\U$ contains all $u\in C^1([t_0,\infty),\Rat{})$, $t_0\ge 0$, which satisfy {\eqref{u:bound}} and the following property:	

\noindent\textbf{Property 1.} For every $t_{0} \ge 0$ there exists a strictly increasing sequence of times $\left\{T_{\nu }\right\}_{\nu\in\mathbb{N}_{0}}$  with $t_0=T_0$ and $\lim_{\nu\to\infty}T_{\nu}=\infty$, such that the following holds. For each $u\in \U(t_0)$, there exists a sequence $\{\hat{t}_{\nu}\}_{\nu\in\mathbb{N}}$ such that  $\hat{t}_{\nu}(=\hat{t}_{\nu}(u))\in(T_{\nu-1},T_{\nu})$ for all $\nu\in \mathbb{N}$,  $\lim_{\nu\to\infty}\hat{t}_{\nu}=\infty$ and
\begin{equation}\label{Property_1}
u(\hat{t}_{\nu } )\ne 0,\forall \nu\in\mathbb{N}.
\end{equation}

\noindent System \eqref{exm:1:system} has the form \eqref{system:triangular} with $f_1(t,x_1,u)=0$, $f_2(t,x_1,x_2,u)=g(t,x_1,x_2,u)-x_2^q$ and $a_{1}(t,x_1,u):=u$ and obviously satisfies hypotheses H1 and H2 with the given input set $\U$ and $M=\Rat{2}$. Moreover, by exploiting our {assumptions \eqref{u:bound} and \eqref{function:g:property},}  it follows that system \eqref{exm:1:system} is $(\Rat{2},\U)$-forward complete. Indeed, by evaluating the time derivative $\dot{V}$ of $V(x_{1},x_{2}):=1/2(x_{1}^{2}+x_{2}^{2})$ along the trajectories of system \eqref{exm:1:system} we get \begin{align*}
\dot{V}(x(t))=&u(t)x_1(t)x_2(t)+x_2(t)g(t,x_1,x_2,u)-x_2^{q+1}(t)\\
\leq&(\bar{u}(t)+\alpha(t)+\beta(t))V(x(t))+\beta(t), \forall t\ge 0,
\end{align*}

\noindent which implies forward completeness. It turns out that system \eqref{exm:1:system}  satisfies the hypotheses H1 and H2, and therefore, according to Proposition \ref{Theorem}, the IDSODP (IDODP) is solvable for \eqref{exm:1:system} with respect to $(M,\U)$.
\end{exm}
\subsection{Comparison with Previous Works}

In \cite{BdTj13a}, systems of the form \eqref{system:triangular} without inputs were considered under the assumption that the system is forward complete and the following implication holds
\begin{align}
\forall t_{0} \in \RgeO ,\, x_{0} \in {\mathbb R}^{n} ,\, i=1,2,\ldots,n-1 \Rightarrow a_{i} (t,y(t))\ne 0\, \, {\rm a.e.\; }t\ge t_{0}. \label{previous_assumption}
\end{align}
One of the main results in \cite{BdTj13a}, establishes that under previous assumptions the \textit{noncausal}-observer design problem is solvable for systems \eqref{system:triangular} by means of a \textit{switching sequence of noncausal observers}.  In addition, in \cite{BdTj13a}, it was also established that if it is a priori known that the initial conditions lie in a bounded subset of $\Rat{n}$, then the noncausal-ODP is solvable for systems \eqref{system:triangular}. A noncausal observer is a system, whose dynamics require sufficient knowledge of future values of both the input and output of the system, in order to estimate its state. As it has been pointed out in \cite{BdTj13a} and \cite{TJ08}, solvability of the ODP by means of a noncausal observer is equivalent to the solvability of the same problem by means of a causal time-delay system.

Obviously, hypothesis H2 is weaker than \eqref{previous_assumption} and therefore Proposition \ref{Theorem} of present paper constitutes a generalization of the main result in \cite{BdTj13a}. In particular, Assumption H2 includes the general case for which the mappings $a_i(\cdot,y(\cdot),u(\cdot))$, $i=1,\ldots,n-1$ may vanish in open subintervals of $\RgeO$.

It should be noted that, extensions of the result in \cite{BdTj13a} and \cite{TJ08} have been established in \cite{TjTd17} for the solvability of the noncausal ODP for systems \eqref{system:triangular}, when a \textit{growth condition} is imposed on the dynamics of the system (the same assumption is imposed in \cite{TJ08}). The main result in \cite{TjTd17} (Proposition 2.1) establishes that the noncausal ODP is solvable for \eqref{system:triangular}, if in addition to  forward completeness of the system, the following conditions hold:

\noindent \textbf{(I)} For each $t_0\in \RgeO$ and $(u,y)\in O(t_0,M)$, where $O(t_0,M)$ is defined in \eqref{output:functions}, an a priori known constant $\xi>0$ exists in such a way that a sequence of times $\{t_{\nu}\}_{\nu\in\mathbb{N}}$ can be determined with
\begin{subequations}
\begin{align}
&\lim_{\nu\to\infty}t_{\nu}=\infty\\
&t_{\nu+1}-t_{\nu}<\xi,\,\,\nu=0,1,2,\ldots
\end{align}
\end{subequations}
such that 
\begin{equation}
a_{i}(t_{\nu},y(t_{\nu}),u(t_{\nu}))\neq0,\,\,\nu=1,2,\ldots,i=1,\ldots,n-1.
\end{equation}
\textbf{(II)} For every $i\in \{1,2,\ldots,n-2\}$ and nonempty open subset $\Delta \subset \RgeO$ for which 
		\begin{subequations}\label{hyp:H2:a_i:zero}
			\begin{equation}\label{hyp:H2:a_i:0}
			a_{i} (t,y(t),u(t))=0,\,\forall t\in \Delta
			\end{equation}
			then
			\begin{equation}\label{hyp:H2:a_ip1:0}
			a_{j} (t,y(t),u(t))=0, \forall j=i+1,i+2,\ldots ,n-1;
			\end{equation}
			moreover, for every $i\ge 2$ and $t\ge 0$ for which \eqref{hyp:H2:a_i:0} holds, we assume:
			\begin{align}
			\max \biggl\{ &\frac{\partial f_{i} }{\partial x_{i} }(t,y(t),\alpha _{i2} ,\ldots ,\alpha _{ii} ,u(t)),
			\ldots ,\frac{\partial f_{n} }{\partial x_{n} } (t,y(t),\alpha _{n2} ,\alpha _{n3} ,\ldots ,\alpha _{nn} ,u(t))\biggr\} \le 0,\nonumber\\
			&\forall \alpha _{j\nu }\in {\mathbb R},\, j=i,\ldots ,n;\, \nu =2,3,\ldots ,j \label{hyp:H2:derivative_max}
			\end{align}
		\end{subequations}

\noindent Obviously (I)  is stronger than H2 and condition (II) is not required for the establishment of Proposition 2.1 of present work. On the other hand, in order to establish Proposition 2.1, we should impose the additional assumption \eqref{u:bound} and to follow a different constructive design than the one adopted in \cite{TjTd17}. In addition, under lack of the growth condition on the dynamics of \eqref{system:triangular:st}, the global state estimation is achieved by means of a switching sequence of observers.

\section{Preliminary Results}\label{section:Observers:general_sys}

The proof of our main result concerning the case \eqref{system:triangular}, is based on some preliminary results concerning the observer design problem for the case of systems \eqref{system:nonlinear} with linear output:
\begin{subequations} \label{system:linearout}
\begin{align}
\dot{x}&=f(t,x,u):=F(t,x,H(t,u)x,u),\nonumber\\
&(t,x,u)\in {\mathbb R}_{\ge 0}\times {\mathbb R}^{n}\times\mathbb{R}^{p}, \label{system:linearout:st} \\
y&=h(t,x,u):=H(t,u)x,\; y\in {\mathbb R}^{k}, \label{system:linearout:out}
\end{align}
\end{subequations}

\noindent where $H:\RgeO\times\mathbb{R}^{p}\to\Rat{k\times n}$ is $C^{0}$ and $F:{\mathbb R}_{\ge 0}\times {\mathbb R}^{n}\times {\mathbb R}^{k}\times\mathbb{R}^{p}\to {\mathbb R}^{n}$ is $C^{0}$ and locally Lipschitz on $x$. We assume that there exist a nonempty subset $M$ of $\Rat{n}$, a continuous function $\bar{u}:\RgeO\to\RgO$, and a set of continuous inputs $\U$ as in \eqref{input:set}, such that for each $t_0\ge 0$, \eqref{u:bound} holds for all $u\in\U(t_0)$, and additionally, system \eqref{system:linearout:st} is $(M,\U)$-forward complete, namely, for each $t_0\ge 0$, $x_0,\in M$ and $u\in\U(t_0,x_0)$, the solution $x(\cdot ):=x(\cdot,t_{0},x_{0};u)$ of \eqref{system:linearout:st} satisfies \eqref{state:bound} for certain $\beta \in NN$. Also, for every $R>0$ and $t\ge 0$ we define:
\begin{equation} \label{output:set}
Y_{R}(t):=\{y\in {\mathbb R}^{k} :y=H(t,u)x,|x|\le\beta(t,R),|u|\le\bar{u}(t)\},
\end{equation}

\noindent where $H(\cdot,\cdot)$ is given in \eqref{system:linearout:out}. It follows that $Y_{R}(t)\ne\emptyset$ for all $t\ge 0$ and that the set-valued map $[0,\infty)\ni t\to Y_{R}(t) \subset {\mathbb R}^{k}$ satisfies the CP. In addition, due to \eqref{state:bound}, for any $t_{0}\ge 0$ and $(u,y)\in O(t_{0},B_{R}\cap M)$ it holds $y(t)\in Y_{R}(t)$ for all $t\ge t_0$, with $O(t_{0},B_{R}\cap M)$ as defined in \eqref{output:functions}. 
For system \eqref{system:linearout} we make the following hypothesis:
\begin{hyp}\label{hypothesis:det}
There exist an integer $\ell\in\mathbb{N}$, a continuous map $A:\RgeO\times\Rat{n}\times\Rat{k}\times\Rat{p}\to\Rat{n\times n}$ and constants $L>1$ and $R>0$, such that the following properties hold:

\noindent \textbf{A1.} For every $\xi>0$ there exists a set-valued map
\begin{equation} \label{map:Q}
[0,\infty) \ni t\to Q_{R}(t):=Q_{R,\xi}(t)\subset {\mathbb R}^{\ell},
\end{equation}

\noindent with $Q_{R}(t)\ne \emptyset$ for all $t\ge 0$, satisfying the CP and such that
\begin{align}
\forall t\ge 0,x,z&\in\Rat{n}\:{\rm with}\:|x|\le\beta (t,R),\: |x-z|\le\xi,\,y\in Y_{R}(t)\:{\rm and}\:u\in \Rat{p}\: {\rm with}\: |u|\le\bar{u}(t)  \,\textrm{it holds}\nonumber \\
 F(t,x,y,u)&-F(t,z,y,u)=A(t,q,y,u)(x-z),\: {\rm for}\: {\rm certain}\: q\in Q_{R}(t), \label{rhs:difference}
\end{align}

\noindent with $Y_{R}(\cdot)$ as given by \eqref{output:set}.

\noindent \textbf{A2.}  For every $\xi>0$, $\tau>0$, $t_0\ge 0$, $\bar{t}_0\ge t_0+\tau$ and $(u,y)\in O(t_0,M)$, there exist a nondecreasing continuous function $\kappa_{R}:=\kappa_{R,\xi,\bar{t}_{0},\tau}\in C^{0}([\bar{t}_{0},\infty);\Rat{})$  satisfying
\begin{equation}\label{map:kappaR}
\lim_{t\to\infty}\kappa_{R}(t)=\infty,
\end{equation}

\noindent a map $Q_{R}:=Q_{R,\xi}$ as in \eqref{map:Q} and a sequence  $\mathcal{A}_{\nu}:=[a_{\nu},b_{\nu}]$,  $\nu\in\mathbb{N}$ of closed intervals with
\begin{equation} \label{closed:intervals:properties}
a_1\ge t_0+\tau,\quad a_{i}<b_{i},\quad a_{i+1}>b_i,\quad b_i-a_i<\tau,\quad \forall i\in\mathbb{N},
\end{equation}

\noindent in such a way that the following hold. There exist a time-varying symmetric matrix $P_{R}:=P_{R,\xi,t_{0},\bar{t}_{0},\tau,y,u}\in C^{1}([\bar{t}_{0},\infty);\Rat{n\times n})$ and a piecewise continuous function $d_{R}:=d_{R,\xi,t_{0},\bar{t}_{0},\tau,y,u}:[\bar{t}_{0},\infty)\to\Rat{}$, both {causal and strongly causal} on each $\mathcal{A}_{\nu}$ with respect to $O(t_{0},M)$, satisfying:
\begin{subequations} \label{A2:hypothesis}
\begin{equation}\label{map:PR:properties}
P_{R}(t)\ge I_{n\times n},\forall t\ge \bar{t}_{0};|P_{R}(\bar{t}_{0})|\le L;
\end{equation}
\begin{equation} \label{map:dR:properties}
\int_{\bar{t}_{0}}^{t} d_{R}(s)ds>\kappa_{R}(t),\forall t\ge \bar{t}_{0};
\end{equation}
\begin{align}
e'P_{R}(t)A(t-\tau,q,y_{\tau}(t),u_{\tau}(t))e+&\tfrac{1}{2}e'\dot{P}_{R}(t)e\le -d_{R}(t)e'P_{R}(t)e, \nonumber \\
\bullet\forall \nu\in\mathbb{N},t\in \mathcal{A}_{\nu},e\in\ker  H(t-\tau,u_{\tau}&(t)),q\in Q_{R}(t-\tau),\; {\rm and} \nonumber \\
\bullet\forall t\in [\bar{t}_{0},\infty)\setminus\cup_{\nu\in\mathbb{N}}\mathcal{A}_{\nu},e\in\Rat{n},q\in &Q_{R}(t-\tau),\;{\rm provided \; that}\;(u,y)\in O(t_{0},B_{R}\cap M ). \label{Lyapunov:ker:inequality}
\end{align}
\end{subequations}
\end{hyp}

The following result, constitutes a modification of Proposition 2.1 in \cite{BdTj13a} and provides results on the state determination of system \eqref{system:linearout}, when it is a priori known that its initial condition lies in the bounded subset $M\cap B_R$ of $\Rat{n}$. 

\begin{prop}\label{Proposition:state_deter}
Consider the system \eqref{system:linearout} and assume that it is $(M,\U)$-forward complete, namely, there exist a nonempty subset $M$ of $\Rat{n}$ and a set of continuous inputs  $\U$ as  in \eqref{input:set}, such that \eqref{state:bound} holds for certain $\beta\in NN$. In addition, assume that Hypothesis \ref{hypothesis:det} is satisfied. Then, the following hold:

\noindent \textbf{(i)} For each $\tau>0$, $t_0\ge 0$, $\bar{t}_0\ge t_0+\tau$ and $(u,y)\in O(t_{0},B_{R}\cap M)$, there exist a piecewise continuous function $\bar{d}_{R}:[\bar{t}_{0},\infty)\to\Rat{}$ satisfying
\begin{subequations}
\begin{equation}\label{map:dRbar:prop}
\bar{d}_{R}(t)\le d_{R}(t),\forall t\ge\bar{t}_{0};\bar{d}_{R}(t)<d_{R}(t),\forall t\in \mathcal{A}_{\nu},\nu\in\mathbb{N},
\end{equation}
\begin{equation}\label{map:dRbar:int}
\int_{\bar{t}_{0}}^{t} \bar{d}_{R}(s)ds\ge\kappa_{R}(t)-1,\forall t\ge \bar{t}_{0};
\end{equation}

\noindent (with $d_{R}(\cdot)$, $\kappa_{R}(\cdot)$  and $\mathcal{A}_{\nu}$ as given in A2) and a {piecewise continuous function $\phi_{R}:[\bar{t}_{0},\infty)\to\RgeO$}, both causal {and strongly causal on each $\mathcal{A}_{\nu}$} with respect to $O(t_0,B_R\cap M)$, such that
\begin{align}
e'P_{R}(t)A(t-\tau,q,y_{\tau}(t),&u_{\tau}(t))e+\tfrac{1}{2}e'\dot{P}_{R}(t)e\le
\phi_{R}(t)|H(t-\tau,u_{\tau}(t))e|^{2}-\bar{d}_{R}(t)e'P_{R}(t)e, \nonumber \\
\forall t\in [\bar{t}_{0},\infty),e&\in\Rat{n},q\in Q(t-\tau).\label{Lyapunov:Ineq:prop}
\end{align}
\end{subequations}
\noindent \textbf{(ii)} Furthermore the following holds. For each $\tau>0$, $t_0\ge 0$, $\bar{t}_0\ge t_0+\tau$, $(u,y)\in O(t_{0},B_{R}\cap M)$ and constant $\xi$ satisfying
\begin{equation}\label{xi}
\xi\ge\sqrt{L}\beta(\bar{t}_{0},R)\exp\left(-\min\{\kappa_{R}(t)-1:t\ge \bar{t}_{0}\}\right),
\end{equation}

\noindent the solution $z(\cdot)$ of system
\begin{subequations}
\begin{align}
\dot{z}(t)=&G_{\bar{t}_{0}}(t,z(t),y_{\tau}(t),u_{\tau}(t))\nonumber\\:=&F(t-\tau,z(t),y_{\tau}(t),u_{\tau}(t))
+\phi_{R}(t)P_{R}^{-1}(t)H'(t-\tau,u_{\tau}(t))(y_{\tau}(t)-H(t-\tau,u_{\tau}(t))z(t))\label{observer:equation} \\
&{\rm with}\: {\rm initial}\: z(\bar{t}_{0})=0 \label{observer:initcond}
\end{align}
\end{subequations}

\noindent with $P_{R}(\cdot)$ as given in A2, is defined for all $t\ge\bar{t}_{0}$. Furthermore, the error $e(t):=x_{\tau}(t)-z(t)$ between the $\tau$ time units delayed value of the trajectory $x(\cdot):=x(\cdot,t_{0},x_{0};u)$ of \eqref{system:linearout:st}, initiated from $x_{0}\in B_{R}\cap M$ at time $t_{0}\ge 0$ and the trajectory $z(\cdot):=z(\cdot,\bar{t}_{0},0;u,y)$ of \eqref{observer:equation} satisfies:
\begin{subequations}
\begin{align}
&|e(t)|<\xi,\forall t\ge\bar{t}_{0}; \label{error:bound} \\
&|e(t)|\le\sqrt{L}\beta(\bar{t}_{0},R)\exp\left(-\kappa_{R}(t)+1\right),\forall t\ge \bar{t}_{0}. \label{error:estimate:2}
\end{align}
\end{subequations}

\noindent It follows from \eqref{map:kappaR} and \eqref{error:estimate:2}, that for $\bar{t}_{0}:=t_{0}+\tau$ the $\tau$-DODP is solvable for \eqref{system:linearout} with respect to $(B_{R}\cap M,\U)$. In particular, the error $e(\cdot)$ between the trajectory $x(\cdot):=x(\cdot,t_{0},x_{0};u)$, $x_{0}\in B_{R}\cap M$, $u\in\U$ of \eqref{system:linearout:st} and the trajectory $z(\cdot):=z(\cdot,\bar{t}_{0},z_{0};u,y)$, $z(\bar{t}_{0})=0$ of the observer {\eqref{observer:equation}} satisfies \eqref{error:conv}. Hence, since $\tau>0$ can be selected arbitrarily small, we also deduce that the IDODP is solvable for \eqref{system:linearout} with respect to $(B_{R}\cap M,\U)$.
\end{prop}

\begin{rem}
Notice, that due to \eqref{map:kappaR} and continuity of $\kappa_R(\cdot)$,  the minimum in \eqref{xi} is {well defined}. In addition, it follows from \eqref{map:dR:properties} that $\kappa_R(\bar{t}_0)<0$ and therefore $\exp(-\min\{\kappa_{R}(t)-1:t\ge \bar{t}_{0}\})>1$.
\end{rem}

The following proposition also constitutes a modification of Proposition 2.2 in \cite{BdTj13a}. It establishes sufficient conditions for the existence of a switching observer (see Definition \ref{dfn:Dsods}) exhibiting the state determination of \eqref{system:linearout}, without any a priori information concerning the initial condition. We make the following hypothesis:
\begin{hyp}\label{hypothesis:switchin}
\textit{There exist a constant $L>1$, an integer $\ell\in{\mathbb N}$ and a map $A:\RgeO\times\Rat{\ell}\times\Rat{k}\times \Rat{p}\to\Rat{n\times n}$, in such a way that for every $R>0$ Hypothesis \ref{hypothesis:det} is fulfilled, namely, both A1 and A2 hold.}
\end{hyp}

\begin{prop}\label{Proposition:switching} In addition to the hypothesis of $(M,\U)$-forward completeness for \eqref{system:linearout:st} {made in Proposition \ref{Proposition:state_deter}}, assume that system \eqref{system:linearout} satisfies Hypothesis \ref{hypothesis:switchin}. Then the IDSODP is solvable for \eqref{system:linearout} with respect to $(M,\U)$. \end{prop}

{The proofs of Propositions \ref{Proposition:state_deter} and \ref{Proposition:switching} are based on a technical result (Lemma \ref{lemma} below) which constitutes a modification of the corresponding result of Lemma 2.1 in \cite{BdTj13a}.  An outline of its proof is provided in the Appendix.}

\noindent Let $k,\ell,n,p,s \in {\mathbb N}$, {$W\subset \Rat{n}$} and consider a pair $(H,A)$ of continuous mappings:
\begin{subequations} \label{maps:HA:aux}
\begin{align}
H:=&H(t,u):\RgeO\times\Rat{p}\to\Rat{s\times n}; \label{map:H:aux} \\
A:=&A(t,q,y,u):\RgeO\times\Rat{\ell}\times\Rat{k}\times\Rat{p}\to\Rat{n\times n}\label{map:A:aux}
\end{align}
\end{subequations}

\noindent Also, consider a set-valued mapping
\begin{equation} \label{map:Q:aux}
[0,\infty)\ni t\to Q(t)\subset {\mathbb R}^{\ell}
\end{equation}

\noindent with $Q(t)\ne \emptyset$ for all $t\ge 0$, satisfying the CP and for each $t_{0}\ge 0$, let $\Omega(t_{0},W)$ be a nonempty set of continuous functions $(u,y):=(u_{t_{0},x_{0}},y_{t_{0},x_{0}}):[t_{0},\infty)\to\Rat{p}\times\Rat{k}$, parameterized by $t_0\in\RgeO$ and $x_{0}\in W$.

\begin{lemma}\label{lemma}
Consider the pair $(H,A)$ of the time-varying mappings in \eqref{maps:HA:aux} and the set valued map $Q(\cdot)$ in \eqref{map:Q:aux}. Also, let $t_0\ge 0$, $\tau>0$, $b>a\ge t_0+\tau$ {with $b-a<\tau$}, $(u,y)\in\Omega(t_0,W)$ and assume that there exist a time-varying symmetric positive definite matrix $P:=P_{t_{0},\bar{t}_{0},\tau,y,u}\in C^{1}([a,b];{\mathbb R}^{n\times n})$ and a function $d:=d_{t_{0},\bar{t}_{0},\tau,y,u}\in C^0([a,b];{\mathbb R})$, both {strongly} causal on $[a,b]$ with respect to $\Omega(t_{0};W)$ and such that
\begin{align}
&e'P(t)A(t-\tau,q,y_{\tau}(t),u_{\tau}(t))e+\tfrac{1}{2}e'\dot{P}(t)e\le -d(t)e'P(t)e, \nonumber \\
&\forall t\in[a,b],e\in\ker H(t-\tau,u_{\tau}(t)),q\in Q(t-\tau) \label{Lyapunov:ker:inequatily:aux}
\end{align}

\noindent and furthermore
\begin{equation} \label{rank:H}
{\rm rank}H(t-\tau,u_{\tau}(t))<n,\forall t\in[a,b].
\end{equation}

\noindent Then for every {strongly} causal on $[a,b]$ with respect to $\Omega(t_{0};W)$ function $\bar{d}:=\bar{d}_{t_{0},\bar{t}_{0},\tau,y,u}\in C^0([a,b];{\mathbb R})$ with
\begin{equation} \label{dbar:lt:d}
\bar{d}(t)<d(t),\forall t\in [a,b],
\end{equation}
\noindent there exists a  {strongly causal on $[a,b]$ with respect to $\Omega(t_{0};W)$} function $\phi:=\phi_{t_{0},\bar{t}_{0},\tau,y,u}\in C^{1}([a,b];{\mathbb R}_{> 0})$ such that
\begin{align}
e'P(t)A(t-\tau,q,&y_{\tau}(t),u_{\tau}(t))e+\tfrac{1}{2}e'\dot{P}(t)e\le \phi(t)|H(t-\tau,u_{\tau}(t))e|^{2}-\bar{d}(t)e'P(t)e, \nonumber \\
\forall t\in[a,b],&e\in\Rat{n},q\in Q(t-\tau). \label{Lyapunov:inequality:aux}
\end{align}
\end{lemma}

We are in a position now to prove Propositions \ref{Proposition:state_deter} and \ref{Proposition:switching} which play a key role in the proof of our main result {(Proposition \ref{Theorem}). For their proof, we will exploit the result of Lemma \ref{lemma}.}

\begin{IEEEproof}[Proof of Proposition \ref{Proposition:state_deter}]
(i)  Let $\tau>0$, ${t}_{0}\ge0$, $ \bar{t}_0\ge t_0+\tau$, $(u,y)\in O(t_0,B_{R}\cap M)$ and a sequence of closed intervals $A_{\nu}=[a_{\nu},b_{\nu}]$, $\nu\in\mathbb{N}$ that satisfy \eqref{closed:intervals:properties}. Also, let $d_{R}(\cdot)$ and $\kappa_{R}(\cdot)$ as given in Property A2 of Hypothesis \ref{hypothesis:det} and define
\begin{equation}\label{map:bar:dR}
\bar{d}_R(t):=\left\{ \begin{array}{ll}
d_{R}(t)-\frac{2}{\pi (1+t^2)}, & t\in \mathcal{A}_{\nu}, \, \forall \nu\in\mathbb{N} \\ 
d_{R}(t), & t\in [\bar{t}_0,\infty)\setminus\cup_{\nu\in\mathbb{N}}\mathcal{A}_{\nu}
\end{array}\right.
\end{equation}
		
\noindent It follows that $\bar{d}_{R}:[\bar{t}_0,\infty)\to\mathbb{R}$ is piecewise continuous, causal and {strongly} causal on each $\mathcal{A}_{\nu}$ with respect to $O(t_0,M)$. In addition, due to \eqref{map:dR:properties}, both \eqref{map:dRbar:prop} and \eqref{map:dRbar:int} are satisfied.

We next construct {a piecewise continuous function $\phi_R:[\bar{t}_0,\infty)\to\RgeO$, being causal and {strongly} causal on each $\mathcal{A}_{\nu}$ with respect to $O(t_0, M)$}, such that \eqref{Lyapunov:Ineq:prop} holds with $\bar{d}_{R}(\cdot)$ as given by \eqref{map:bar:dR}. {First}, we successively apply the result of Lemma \ref{lemma} to determine the {restriction} of $\phi_{R}(\cdot)$ on {each closed interval} $\mathcal{A}_{\nu}$, $\nu\in\mathbb{N}$. {In particular}, by taking into account \eqref{closed:intervals:properties}, \eqref{Lyapunov:ker:inequality}, the second inequality of \eqref{map:dRbar:prop} and \eqref{map:bar:dR} it follows that for each $\nu\in\mathbb{N}$ the hypotheses of Lemma \ref{lemma} are satisfied with  $W:=B_{R}\cap M$ and $[a,b]:=\mathcal{A}_{\nu}$. Therefore, for each $\nu\in \mathbb{N}$, there exists a function $\phi_{R,\nu}\in C^{1}(\mathcal{A}_{\nu};\mathbb{R}_{>0})$, being  {strongly} causal on each $\mathcal{A}_{\nu}$ with respect to $O(t_0;M)$ in such a way that	
\begin{align}
e'P_R(t)A(t-&\tau,q,y_{\tau}(t),u_{\tau}(t))e+\tfrac{1}{2}e'\dot{P}_R(t)e\le\phi_{R,\nu}(t)|H(t-\tau,u_{\tau}(t))e|^{2}-\bar{d}_R(t)e'P_R(t)e, \nonumber \\
\forall t\in \mathcal{A}_{\nu},&\nu\in\mathbb{N},e\in\Rat{n},q\in Q(t-\tau). \label{Lyapunov:inequality:prop:1a}
\end{align}

\noindent {Then, we select a non-negative piecewise continuous function} $\phi_{R}:[\bar{t}_0,\infty)\to\Rat{}$, such that $\phi_{R}(t)=\phi_{R,\nu}(t)$ for all $t\in\mathcal{A}_{\nu}$, $\nu\in\mathbb{N}$. Thus, by virtue of \eqref{Lyapunov:ker:inequality}, \eqref{map:dRbar:prop}, \eqref{map:bar:dR}, \eqref{Lyapunov:inequality:prop:1a}, the fact that $\phi_R(t)\geq 0$ for all $t\ge \bar{t}_0$, and by recalling that $(u,y)\in O(t_{0},B_{R}\cap M)$, it follows that
\begin{align}
e'P_{R}(t)A(t-\tau,&q,y_{\tau}(t),u_{\tau}(t))e+\tfrac{1}{2}e'\dot{P}_{R}(t)e\leq \phi_{R}(t)|H(t-\tau,u_{\tau}(t))e|^{2}-\bar{d}_{R}(t)e'P_{R}(t)e, \nonumber \\
\forall t\in [\bar{t}_{0},\infty),e & \in\Rat{n},{q\in Q(t-\tau).}\label{Lyapunov:Ineq:prop:1b}
\end{align}

\noindent This completes the proof of the first statement of Proposition \ref{Proposition:state_deter}.

\noindent In order to prove part (ii), 
let $\tau>0$, ${t}_{0}\ge0$, $ \bar{t}_0\ge t_0+\tau$, $u\in \U$, $(u,y)\in O(t_0,B_R\cap M)$ and {a constant  $\xi$ as in} \eqref{xi}. According to Hypothesis \ref{hypothesis:det}, there exist a sequence $\{\mathcal{A}_{\nu}\}_{\nu\in\mathbb{N}}$ of {closed intervals satisfying \eqref{closed:intervals:properties}}, a time-varying symmetric matrix $P_{R}\in C^{1}([\bar{t}_0,\infty);\mathbb{R}^{n\times n})$ and a piecewise continuous function $d_{R}:[\bar{t}_0,\infty)\to\mathbb{R}$, both {strongly} causal on each  $\mathcal{A}_{\nu}$, $\nu\in\mathbb{N}$, with respect to $O(t_0,B_R\cap M)$, {and such that \eqref{map:PR:properties}-\eqref{Lyapunov:ker:inequality} hold}. Then, statement (i) of the proposition asserts that there exist piecewise continuous {functions} $\bar{d}_{R}:[\bar{t}_0,\infty)\to\mathbb{R}$ and ${\phi_{R}:[\bar{t}_0,\infty)\to \mathbb{R}_{\ge 0}}$, both causal {and {strongly} causal on each $\mathcal{A}_{\nu}$} with respect to $O(t_0,B_{R}\cap M)$, such that \eqref{map:dRbar:prop}, \eqref{map:dRbar:int} and \eqref{Lyapunov:Ineq:prop} hold. By invoking \eqref{system:linearout:out}, the error equation $e(t)=x_{\tau}(t)-z(t)$ between \eqref{system:linearout:st} and \eqref{observer:equation} is written:
\begin{align}
\dot{e}(t)=& F(t-\tau,x_{\tau}(t),y_{\tau}(t),u_{\tau}(t))-F(t-\tau,z(t),y_{\tau}(t),u_{\tau}(t))-\phi_R(t)P_{R}^{-1}(t)\nonumber\\
& \times H'(t-\tau,u_{\tau}(t))(y_{\tau}(t)-H(t-\tau,u_{\tau}(t))z(t))\nonumber\\
 =& F(t-\tau,x_{\tau}(t),y_{\tau}(t),u_{\tau}(t))-F(t-\tau,z(t),y_{\tau}(t),u_{\tau}(t))-\phi_R(t)P_{R}^{-1}(t)\nonumber\\
& \times  H'(t-\tau,u_{\tau}(t))(H(t-\tau,u_{\tau}(t))x-H(t-\tau,u_{\tau}(t))z(t))\nonumber\\
 =& F(t-\tau,x_{\tau}(t),y_{\tau}(t),u_{\tau}(t))-F(t-\tau,z(t),y_{\tau}(t),u_{\tau}(t))-\phi_{R}(t)P_{R}^{-1}(t)\nonumber\\
& \times H'(t-\tau,u_{\tau}(t))H(t-\tau,u_{\tau}(t))e, \, \, t\ge \bar{t}_0 \label{error:derivative}
\end{align}

\noindent Next, recall \eqref{observer:initcond}, namely, that $z(\bar{t}_0)=0$ and let $[\bar{t}_0,T_{\max})$ be the maximum right-interval of existence of the solution $e(\cdot)$ of \eqref{error:derivative} with initial condition
\begin{equation}\label{error:initcond}
e(\bar{t}_0)=x(\bar{t}_0-\tau,t_0,x_0;u)-z(\bar{t}_0)=x(\bar{t}_0-\tau,t_0,x_0{;u}).
\end{equation}
	
\noindent Then, from \eqref{state:bound}, \eqref{xi}, \eqref{error:initcond}, the fact that $L>1$ (see Hypothesis \ref{hypothesis:det}), $\beta\in NN$, $x_0\in B_R\cap M$ {and by virtue of Remark 3.2}, it follows that
\begin{align}
|e(\bar{t}_0)|&=|x(\bar{t}_0-\tau,t_0,x_0;u)|\leq {\beta(\bar{t}_0-\tau,|x_0|)\le\beta(\bar{t}_0,R)}\nonumber\\
&< {\beta(\bar{t}_0,R)\sqrt{L}\exp(-\min\{\kappa_{R}(t)-1:t\ge \bar{t}_0\})}\leq\xi. \label{error:ineq:t0}
\end{align}

\noindent We claim that $|e(t)|<\xi$ for every $t\in[\bar{t}_0,T_{\max})$ and therefore $T_{\max}=\infty$. Indeed, suppose on the contrary that there exists a time $\hat{t}\in[\bar{t}_0,T_{\max})$ such that
\begin{subequations}\label{Error:contrad}
\begin{align}
|e(\hat{t})|&=\xi\\
|e(t)|&<\xi, \, \, \forall t\in [\bar{t}_0,\hat{t})
\end{align}
\end{subequations}

\noindent By recalling \eqref{output:set},  and taking into account \eqref{Error:contrad} and the fact that $e(t)=x_{\tau}(t)-z(t)$, {$y_{\tau}(t)\in Y_R(t-\tau)$ and $|u_{\tau}(t)|\le \bar{u}(t-\tau)$, it follows from \eqref{rhs:difference} of Property A1 applied with $t:=t-\tau$}, that for each $t\in [\bar{t}_0,\hat{t}]$ there exists $q\in Q_R(t-\tau)$ such that
\begin{align}
F(t-\tau&,x_{\tau}(t),y_{\tau}(t),u_{\tau}(t))-F(t-\tau,z(t),y_{\tau}(t),u_{\tau}(t))=A(t-\tau,q,y_{\tau}(t),u_{\tau}(t))({x_{\tau}(t)}-z(t)).\label{rhs:diff:prop}
\end{align}

\noindent Thus, by evaluating the time derivative $\dot{V}$ of $V(t,e):=\frac{1}{2}e'P_{R}(t)e$, $e\in \mathbb{R}^n$ along the trajectories $e(\cdot)=e(\cdot,\bar{t}_0,e(\bar{t}_0);{x(\cdot),u(\cdot)})$ of \eqref{error:derivative} and by exploiting \eqref{rhs:diff:prop} we have
\begin{align}
\dot{V}(t,e(t))=&\frac{1}{2}\frac{d}{dt}(e'(t)P_{R}(t)e(t))\nonumber\\
=&\frac{1}{2}\dot{e}'(t)P_{R}(t)e(t)+\frac{1}{2}e'(t)\dot{P}_{R}(t)e(t)+\frac{1}{2}e'(t)P_{R}(t)\dot{e}(t)\nonumber\\
=&\frac{1}{2}e'(t)\dot{P}_{R}(t)e(t)+e'(t)P_{R}(t)(F(t-\tau,x_{\tau}(t),y_{\tau}(t),u_{\tau}(t))-F(t-\tau,z(t),y_{\tau}(t),u_{\tau}(t)))\nonumber\\
&- e'(t)P_R(t)\phi_R(t)P^{-1}_R(t)H'(t-\tau,u_{\tau}(t))H(t-\tau,u_{\tau}(t))e(t)\nonumber\\
=&\frac{1}{2}e'(t)\dot{P}_{R}(t)e(t)+e'(t)P_{R}(t)A(t-\tau,q,y_{\tau}(t),u_{\tau}(t))e(t)\nonumber\\
&- \phi_{R}(t)|H(t-\tau,u_{\tau}(t))e(t)|^2,\forall t\in [\bar{t}_0,\hat{t}], q\in Q_R(t-\tau). \label{LyapErr:derivative}
\end{align}

\noindent Therefore, by \eqref{Lyapunov:Ineq:prop} and \eqref{LyapErr:derivative} it follows that
\begin{equation}\label{LyapErr:ineq}
\dot{V}(t,e(t))\leq -2\bar{d}_R(t)V(t,e(t)), \, \, \forall t\in [\bar{t}_{0},\hat{t}].
\end{equation}

\noindent From \eqref{map:dRbar:prop}, \eqref{LyapErr:ineq} and the definition of $V(\cdot,\cdot)$ above we find that
	\begin{equation}\label{LyapErr:ineq2}
	V(t,e(t))=\frac{1}{2}e'(t)P_R(t)e(t)\leq \frac{1}{2}e'(\bar{t}_0)P_{R}(\bar{t}_0)e(\bar{t}_0)\exp\Bigl(-\int_{\bar{t}_0}^{t}2\bar{d}_R(s)ds\Bigr), \forall t\in[\bar{t}_0,\hat{t}].
	\end{equation}

\noindent {Thus, it follows from the first inequality of} \eqref{map:PR:properties} and \eqref{LyapErr:ineq2} that
\begin{align}
e'(t)P_{R}(t)e(t) & \leq e'(\bar{t}_0)P_{R}(\bar{t}_0)e(\bar{t}_0)\exp\Biggl(-\int_{\bar{t}_0}^{t}2\bar{d}_R(s)ds\Biggr)\Longrightarrow \nonumber\\
|e(t)| & \leq |e(\bar{t}_0)|\sqrt{|P_{R}(\bar{t}_0)|}\exp\Biggl(-\int_{\bar{t}_0}^{t}\bar{d}_R(s)ds\Biggr),\forall t\in [\bar{t}_0,\hat{t}].\label{Error:ineq}
\end{align}

\noindent In addition, from \eqref{state:bound} {and} \eqref{error:initcond} we have that
\begin{equation}\label{Error:bound}
|e(\bar{t}_0)|=|x(\bar{t}_0-\tau,t_0,x_0;u)|\leq \beta(\bar{t}_0,R).
\end{equation}

\noindent Consequently, by taking into account {the second inequality of} \eqref{map:PR:properties}, \eqref{map:dRbar:int}, \eqref{xi}, \eqref{error:ineq:t0}, \eqref{Error:ineq} and \eqref{Error:bound} we get that
\begin{align}
|e(t)|\leq&\beta(\bar{t}_0,R)\sqrt{L}\exp{\Biggl(-\int_{\bar{t}_0}^{t}\bar{d}_R(s)ds\Biggr)} \nonumber\\
\leq&\beta(\bar{t}_0,R)\sqrt{L}{\exp(-\kappa_{R}(t)+1)} \nonumber\\
<&\beta(\bar{t}_0,R)\sqrt{L}\exp(-\min\{\kappa_{R}(t)-1: t\in [\bar{t}_0,\hat{t}]\})\leq \xi, \forall t\in [\bar{t}_0,\hat{t}].
\end{align}

\noindent Thus $|e(\hat{t})|<\xi$, which contradicts \eqref{Error:contrad}. Therefore the solution $e(\cdot)=e(\cdot,\bar{t}_0,e(\bar{t}_0);x(\cdot),u(\cdot))$  of \eqref{error:derivative} satisfies \eqref{error:bound}, namely $|e(t)|<\xi$, for every $t\in [\bar{t}_0,T_{\max})$, hence $T_{\max}=\infty$. {The latter implies that \eqref{Error:ineq} holds for all $t\ge t_0$. Hence, by recalling again the second inequality of \eqref{map:PR:properties}, \eqref{map:dRbar:int}, \eqref{xi}, \eqref{Error:ineq} and \eqref{Error:bound}, it follows that \eqref{error:estimate:2} is fulfilled as well.} Finally, by taking into account \eqref{map:kappaR}, \eqref{map:dRbar:int} and \eqref{error:estimate:2} we conclude that $\lim_{t\rightarrow \infty}e(t)=0$ and therefore the $\tau$-DODP is solvable for \eqref{system:linearout} with respect to $(B_{R}\cap M,\U)$.
\end{IEEEproof}

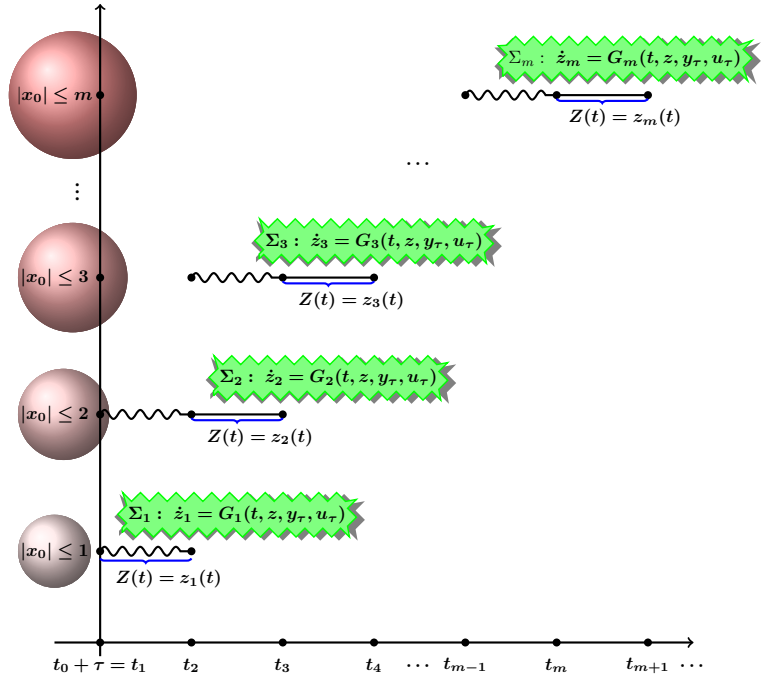
\begin{figure}[!t]
\centering\resizebox{!}{9cm}{%
\begin{tikzpicture}[thick,scale=0.9]
\tikzstyle{blocknew} = [rectangle, circular drop shadow={shadow scale=1.05}, decorate, decoration=zigzag, draw=green, thick, fill=green!50, text centered, minimum height=2em]
\tikzstyle{textblock} = [rectangle, thick, text width=11cm, text centered, minimum height=2em]

\shade[ball color=red!10] (-1,2) circle (0.8cm);
\shade[ball color=red!20] (-0.8,5) circle (1.cm);
\shade[ball color=red!30] (-0.6,8) circle (1.2cm);
\shade[ball color=red!40] (-0.6,12) circle (1.4cm);

\draw[->,very thick] (-1,0) to (13,0);
\draw[->,very thick] (0,-0.3) to (0,14);

\fill (0,0) circle (2.5pt);
\fill (2,0) circle (2.5pt);
\fill (4,0) circle (2.5pt);
\fill (6,0) circle (2.5pt);
\fill (8,0) circle (2.5pt);
\fill (10,0) circle (2.5pt);
\fill (12,0) circle (2.5pt);

\fill (0,2) circle (2.5pt);
\fill (2,2) circle (2.5pt);

\fill (0,5) circle (2.5pt);
\fill (2,5) circle (2.5pt);
\fill (4,5) circle (2.5pt);

\fill (0,8) circle (2.5pt);
\fill (2,8) circle (2.5pt);
\fill (4,8) circle (2.5pt);
\fill (6,8) circle (2.5pt);

\fill (0,12) circle (2.5pt);
\fill (8,12) circle (2.5pt);
\fill (10,12) circle (2.5pt);
\fill (12,12) circle (2.5pt);

\draw[decorate,decoration={coil,aspect=0}, very thick] (0,2) -- (2,2);
\draw[decorate,decoration={brace,mirror,raise=4pt},blue, very thick] (0,2) -- (2,2);

\draw[decorate,decoration={coil,aspect=0}, very thick] (0,5) -- (2,5);
\draw[very thick] (2,5) to (4,5);
\draw[decorate,decoration={brace,mirror,raise=2pt},blue, very thick] (2,5) -- (4,5);

\draw[decorate,decoration={coil,aspect=0}, very thick] (2,8) -- (4,8);
\draw[very thick] (4,8) to (6,8);
\draw[decorate,decoration={brace,mirror,raise=2pt},blue, very thick] (4,8) -- (6,8);

\draw[decorate,decoration={coil,aspect=0}, very thick] (8,12) -- (10,12);
\draw[very thick] (10,12) to (12,12);
\draw[decorate,decoration={brace,mirror,raise=2pt},blue, very thick] (10,12) -- (12,12);

\node at (-1, 2) {$\bm{|x_{0}|\le 1}$};
\node at (-1, 5) {$\bm{|x_{0}|\le 2}$};
\node at (-1, 8) {$\bm{|x_{0}|\le 3}$};
\node at (-1, 12) {$|\bm{x_{0}|\le m}$};

\node at (0, -0.5) {$\bm{t_{0}+\tau=t_{1}}$};
\node at (2, -0.5) {$\bm{t_{2}}$};
\node at (4, -0.5) {$\bm{t_{3}}$};
\node at (6, -0.5) {$\bm{t_{4}}$};
\node at (7, -0.5) {$\bm{\cdots}$};
\node at (8, -0.5) {$\bm{t_{m-1}}$};
\node at (10, -0.5) {$\bm{t_{m}}$};
\node at (12, -0.5) {$\bm{t_{m+1}}$};
\node at (13, -0.5) {$\bm{\cdots}$};

\node at (7, 10.5) {$\bm{\cdots}$};
\node at (-0.5, 10) {$\bm{\vdots}$};

\node at (1.5,1.4 ) {$\bm{Z(t)=z_{1}(t)}$};
\node at (3,2.8) [blocknew] {$\bm{\Sigma_1:\,\,\dot{z}_{1}=G_{1}(t,z,y_{\tau},u_{\tau})}$};

\node at (3.5,4.5 ) {$\bm{Z(t)=z_{2}(t)}$};
\node at (5,5.8) [blocknew] {$\bm{\Sigma_2:\,\,\dot{z}_{2}=G_{2}(t,z,y_{\tau},u_{\tau})}$};

\node at (5.5,7.5 ) {$\bm{Z(t)=z_{3}(t)}$};
\node at (6,8.8) [blocknew] {$\bm{\Sigma_3:\,\,\dot{z}_{3}=G_{3}(t,z,y_{\tau},u_{\tau})}$};

\node at (11.5,11.5 ) {$\bm{Z(t)=z_{m}(t)}$};
\node at (11.5,12.8) [blocknew] {$\Sigma_m:\,\,\bm{\dot{z}_{m}=G_{m}(t,z,y_{\tau},u_{\tau})}$};
\end{tikzpicture} }
\caption{Switching Sequence of Observers.} \label{fig:switching}
\end{figure}

We proceed to the proof of Proposition \ref{Proposition:switching} which establishes solvability of the ODP for systems \eqref{system:nonlinear} by means of a delayed switching sequence of observers. Its proof is similar to the proof of Proposition 2.2 in \cite{BdTj13a}. In particular, under A1 and A2, we can construct the desired observers according to Definition \ref{dfn:Dsods}. The main idea behind the proof is illustrated through Fig. \ref{fig:switching} which depicts a strictly increasing sequence of times, with $t_0$ being the initial time of the system and a switching sequence of observers. Specifically, we successively apply the estimates of Proposition \ref{Proposition:state_deter} with $R=1,2,3,\ldots$ {by pretending} that the initial condition of the system belongs to the intersection of $M$ and the ball $B_R$ in $\Rat{n}$ with center 0 and radius $R=1,2,3,\ldots$ respectively. We next focus our informal discussion on the design of observer $\Sigma_2$. By applying Proposition \ref{Proposition:state_deter} with  $R=2$, we may consider the observer $\Sigma_1$ whose existence is guaranteed by Proposition \ref{Proposition:state_deter},  which can perform the asymptotic estimation of the state, provided that the initial state $x_0$ belongs to the ball $B_2$ with center 0 and radius $R=2$. Despite the fact that $x_0$ does not necessarily belong to that ball, we select $t_2$ so that the error is sufficiently small in the case where indeed $x_0\in B_2$. At time $t_2$ we start the design of observer $\Sigma_3$ under the same reasoning for the ball $B_3$ of center 0 and radius $R=3$, and choose time $t_3$ respectively, as if it were the case that $x_0\in B_3$. At $t_3$ we terminate operation of the observer $\Sigma_2$ and proceed with the design of $\Sigma_4$. Proceeding with this recursive design, we can achieve the global state estimation for each initial condition $x_0\in M\cap \Rat{n}$ through the piecewise continuous function  $Z(t):=z_m(t)$, $t\in[t_m,t_m+1)$, $m=1,2,\ldots$, where $z_{m}(\cdot)$ is the trajectory of observer $\Sigma_m$. We also emphasize that the times $t_m$ are chosen in such a way that the estimates of Proposition \ref{Proposition:state_deter} guarantee that  $Z$ converges to the solution of the system.

\begin{IEEEproof}[Proof of Proposition \ref{Proposition:switching}]
{Without any loss of generality, we may assume that $0\in cl M$ and therefore $B_R\cap M\neq \emptyset$ for every $R>0$, which according to Hypothesis \ref{hypothesis:switchin}, implies that A1 and A2 of Hypothesis \ref{hypothesis:det}  hold for every $R=1,2,\ldots$.  Consider the system \eqref{system:linearout}} initiated at $t_0\ge 0$ and let $\tau>0$, $(u,y)\in O(t_0,M)$, $L>1$ and $R>0$. We proceed to the construction of a sequence of times $(t_{m})_{m\in\mathbb{N}}$ and a sequence of continuous and causal with respect to $O(t_0,M)$ mappings $(G_{m})_{m\in\mathbb{N}}$ that satisfy the requirements of  Definition \ref{dfn:Dsods} for the solvability of the IDSODP. The proof is carried out by induction and is based on the following claim.
	
\noindent\textbf{Claim 1(Induction Hypothesis):} Given $L>1$, $t_{0}\ge0$, $\tau>0$ and $(u,y)\in O(t_0,M)$ as above, for any $m\in {\mathbb N}$ there exist:
	
\noindent$\bullet$ positive constants  $\xi _{m} $ and ${t_{m} :=t_{m,t_{0} ,t_{m-1} ,\tau ,\xi _{m}}}$, in such a way that the sequence $(t_{m} )_{m\in {\mathbb N}}$ satisfies:
\begin{subequations}\label{time:prop}
\begin{align}
t_{m+1} & \ge t_{m} +1, m=1,2,\ldots  \label{ineq:req:time} \\
\mathop{\lim }\limits_{m\to \infty } t_{m}&  =\infty, \label{req:time:infty}
\end{align}
\end{subequations}

\noindent with $t_1=t_0+\tau$.
	
\noindent$\bullet $ a set-valued mapping $Q_{m}(\cdot):=Q_{m,t_{0},\xi _m}(\cdot)$ satisfying the CP {(see Notations)} such that \eqref{rhs:difference} is fulfilled with $R:=m$, $t_{0}\ge 0$ and $\xi :=\xi _{m} $;
	
\noindent$\bullet $ a nondecreasing continuous function $\kappa_m:[{t}_{0}+\tau,\infty)\to \Rat{}$ satisfying \eqref{map:kappaR}, namely $\lim_{t\to\infty}\kappa_{m}(t)=\infty$, a symmetric matrix $P_{m} :=P_{m,t_{0},t_{m-1},\tau ,\xi _{m} ,y,u} \in C^{1} ([t_{m-1} ,\infty );{\mathbb R}^{n\times n} )$, {piecewise continuous functions $\bar{d}_{m} :=\bar{d}_{m,t_{0} ,t_{m-1} ,\tau ,\xi _{m} ,y,u}: [t_{m-1} ,\infty )\to{\mathbb R}$ and $\phi _{m} :=\phi _{m,t_{0} ,t_{m-1} ,\tau ,\xi _{m} ,y,u}:[t_{m-1} ,\infty )\to{\mathbb R}_{>0} $}, all causal with respect to $O(t_0,M)$ and in such a way that the following hold:
\begin{subequations}\label{switching:properties}
\begin{equation}\label{map:s:Pm}
P_{m} (t)\ge I_{n\times n}, \forall t\ge t_{m-1} ;|P_{m} (t_{m-1} )| \le L,
\end{equation}
\begin{equation}\label{int:s:dm}
\int_{t_{m-1} }^{t }\bar{d}_{m} (s)ds \geq {\kappa_m}(t)-1, \forall t\geq t_{m-1},
\end{equation}
\begin{equation}\label{Lyap:s:Ineq}
\begin{array}{l} e'P_{m} (t)A(t-\tau,q,y_{\tau}(t),u_{\tau}(t))e+{\tfrac{1}{2}} e'\dot{P}_{m} (t)e-\phi _{m} (t)|H(t-\tau,u_{\tau}(t))e|^{2} \le -\bar{d}_{m} (t)e'P_{m} (t)e, \\ \, \, \forall \, t\ge t_{m-1} ,\, e\in {\mathbb R}^{n} ,\, q\in Q_{m} (t-\tau), {\rm provided}\, {\rm that}\, {(u,y)\in {O(t_0,B_{m}\cap M).}} \end{array}
\end{equation}
\end{subequations}

\noindent Specifically, for each $m\in {\mathbb N}$, the desired constants $\xi _{m} $ and  $t_{m} $ are defined as follows:
\begin{align}
\xi _{i} :=&{\beta(t_{i-1},i)}\sqrt{L} {\exp (-\min\{\kappa_{i}(t)-1:t\in[t_{i-1},t_i]\})}, i=1,\ldots,m,\label{xi:i}
\end{align}
\begin{align}
t_{i} :=\min \{& t\ge t_{i-1} +1:\exp \left(-\kappa_{i}(\bar{t})+1 \right)\le \frac{1}{i \beta(t_{i-1},i)\sqrt{L}},\forall \bar{t}\ge t \},  \nonumber\\
&i=2,\ldots,m\; {\rm for}\; m\ge 2;t_{1} :=t_{0}+\tau. \label{recursive:time}
\end{align}
	
\noindent \textbf{Establishment of Claim 1 {for $m=1$:}}  Due to Hypothesis \ref{hypothesis:switchin}, we apply Properties A1 and A2 with $R=m:=1$. These properties assert that if we define $t_{1} :=t_{0}+\tau $ and $\xi _{1} :=\beta (t_{0}+\tau,1)\sqrt{L} {\exp (-\min\{\kappa_{1}(t) }-1:t\in[t_0,t_1]\})$, then there exists a set-valued map $Q_{\, 1} :=Q_{\, 1,t_{0} ,\xi _{1} } $ satisfying the CP in such a way that \eqref{rhs:difference} holds. In addition, for $\bar{t}_0:=t_1(=t_0+\tau)$, a piecewise continuous function $d_{1}:[\bar{t}_{0} ,\infty )\to{\mathbb R}$ and a time varying symmetric matrix $P_{1} \in C^{1} ([\bar{t}_{0} ,\infty );{\mathbb R}^{n\times n} )$ can be found, both causal with respect to $O(t_0,M)$, such that for $R=1$ conditions \eqref{map:PR:properties}, \eqref{map:dR:properties} hold, and also \eqref{Lyapunov:ker:inequality} is satisfied, provided that $(u,y)\in {O(t_0,B_1\cap M)}$. Therefore, for $R=1$, Proposition \ref{Proposition:state_deter} asserts that there exist {piecewise continuous functions $\bar{d}_{1} :[\bar{t}_{0} ,\infty )\to\Rat{}$ and  $\phi _{1}:[\bar{t}_{0} ,\infty )\to{\mathbb R}_{\ge 0} $} both causal with respect to $O(t_{0},M)$, which satisfy {\eqref{map:s:Pm}}, \eqref{int:s:dm} and \eqref{Lyap:s:Ineq} with $m:=1$.
	
\noindent \textbf{General Induction Step:} Suppose that \eqref{map:s:Pm}-\eqref{Lyap:s:Ineq} are fulfilled for all $i=1,\ldots,m$ for certain integer $m\in {\mathbb N}$ and for appropriate $t_{i} $ satisfying, $t_{i+1} \ge t_{i} +1$, $i=1,2,\ldots,m-1$. We apply again Properties A1 and A2 with $R:=m+1$. If we define
\begin{align*}
\xi _{m+1} :=\beta (t_{m},m+1)\sqrt{L} &\exp(-\min\{\kappa_{m+1}(t)-1:t\in[t_m,t_{m+1}]\}),
\end{align*}

\noindent then a set-valued map $Q_{m+1} :=Q_{m+1,t_{0} ,\xi _{m+1} } $ satisfying the CP can be found, such that \eqref{rhs:difference} is satisfied. In addition,  for $\bar{t}_{0} :=t_{m} \ge t_{0}+\tau $, conditions \eqref{map:PR:properties}, \eqref{map:dR:properties} and \eqref{Lyapunov:ker:inequality} hold for $R:=m+1$, provided that $(u,y)\in {O(t_0,B_{m+1}\cap M)}$. By invoking again Proposition \ref{Proposition:state_deter}, there exist {piecewise continuous functions $\bar{d}_{m+1}:[t_{m} ,\infty )\to{\mathbb R}$ and $\phi _{m+1}:[t_{m} ,\infty )\to{\mathbb R}_{\ge 0} $} both causal with respect to $O(t_{0}, M)$, satisfying \eqref{map:s:Pm}, \eqref{int:s:dm} and \eqref{Lyap:s:Ineq} with $m:=m+1$. Finally, the desired $t_{m+1} $ is given by \eqref{recursive:time} with $i:=m+1$ , namely:
\begin{align}
t_{m+1} :=&\min \{t\ge t_{m} +1:\exp \left(-{\kappa_{m+1}}(t)+1  \right)\le \frac{1}{(m+1)\beta(t_{m},m+1)\sqrt{L} }\}.\label{recursive:time:m+1}
\end{align}

\noindent This completes the establishment of our claim.
	
\noindent \textbf{Design of the switching observer:} {By employing Claim 1 we explicitly construct the desired switching observer} exhibiting the state determination of \eqref{error:conv}, according to the Definition \ref{dfn:Dsods}. {The procedure is similar to the methodology employed in \cite{BdTj13a}.} Consider for each $m\in {\mathbb N}$ the system\begin{subequations}\label{switching:observer}
	\begin{equation}\label{switching:obs}
	\Sigma_m:\,\,\,	\dot{z}_{m}(t) =G_{m} (t,z_{m}(t) ,y_{\tau}(t),u_{\tau}(t)), t\in [t_{m-1} ,t_{m+1} ]
	\end{equation}
with initial
	\begin{equation}\label{switching:obs:initial}
		z_{m} (t_{m-1} )=0\end{equation}
	\begin{align}\label{switching:obs:def}
		G_{m} (t,z_{m},y,u)\,\left\{\begin{array}{l} :=F(t-\tau,z,y,u)+\phi _{m} (t)P_{m}^{-1} (t)H'(t-\tau,u)(y-H(t-\tau,u)z),\\\,\,\,\,\,\,\, \, \, \, \, {\rm for}\, \, |z|\, \le \zeta _{m} ,\, \, t\in [t_{m-1} ,t_{m+1} ],\\{} \\ :={(F(t-\tau,z,y,u)}+\phi _{m} (t)P_{m}^{-1} (t)H'(t-\tau,u)(y-H(t-\tau,u)z))\frac{2\zeta _{m} -|z|}{\zeta _{m} }  \\ {\, \, \, \, \, \, \, \, \, \, \, \, \, \, \, \, \, \, \, \, {\rm for}\, \, \zeta _{m} \le \, |z|\, \le 2\zeta _{m} ,\, \, t\in [t_{m-1} ,t_{m+1} ],} \\ {} \\ {:=0,\, \, \, \, \, \, \, \, \, \, \, \, {\rm for}\, \, |z|\, \ge 2\zeta _{m} ,\, \, t\in [t_{m-1} ,t_{m+1} ]} \end{array}\right.
	\end{align}
	\begin{equation}\label{switching:obs:zeta}
		\zeta _{m} :=\beta (t_{m+1} ,m)+\xi _{m}
	\end{equation}
\end{subequations}

\noindent where the sequences $(t_{m} )_{m\in {\mathbb N}} $, $(\xi _{m} )_{m\in {\mathbb N}} $, $(\phi _{m} )_{m\in {\mathbb N}} $ and $(P_{m} )_{m\in {\mathbb N}} $ are determined in Claim 1. Notice that $G_{m} (\cdot )$ is $C^{0} ([t_{m-1} ,t_{m+1} ]\times {\mathbb R}^{n} \times {\mathbb R}^{k}\times\Rat{p} ;{\mathbb R}^{n} )$ and further it is locally Lipschitz on $z\in {\mathbb R}^{n}$. Moreover, $G_{m} :[t_{m-1} ,t_{m+1} ]\times {\mathbb R}^{n} \times {\mathbb R}^{k}\times\Rat{p} \to {\mathbb R}^{n} $ is bounded, therefore for any initial $z_{m} (t_{0} )\in {\mathbb R}^{n} $, $t_{0} \in [t_{m-1} ,t_{m+1} ]$ the corresponding solution of \eqref{switching:obs} is defined for all $t\in [t_{m-1} ,t_{m+1} ]$.

\noindent\textbf{Establishment of \eqref{error:conv}:} Let $Z:[t_{0} ,\infty )\to {\mathbb R}^{n} $ as defined {in Definition \ref{dfn:Dsods}}, namely, $Z(t):=z_{m} (t)$, $t\in [t_{m} ,t_{m+1} )$, $m\in {\mathbb N}$, where for each $m\in {\mathbb N}$ the map $z_{m} (\cdot )$ is the solution of \eqref{switching:observer}. Notice that for any initial state $x_{0} \in M $ of \eqref{system:linearout:st}, there exists $m_{0} \in {\mathbb N}$ with $m_{0} \ge 2$ such that
\begin{equation} \label{bound:m0}
m_{0} \ge \, |x_{0} |.
\end{equation}
	
\noindent Let $m\ge m_{0} $ and notice that, due to \eqref{switching:obs:initial}, there exists a time $\bar{t}\in (t_{m-1} ,t_{m+1} ]$ such that $|z_{m} (t)|\, <\zeta _{m} $ for all $t\in [t_{m-1} ,\bar{t})$.
	
\noindent \textbf{Claim 2}:

\noindent We claim that
\begin{equation}\label{claim2:prop2}
|z_{m} (t)|<\zeta _{m} , \forall t\in [t_{m-1} ,t_{m+1} ),
\end{equation}
namely, it holds $\bar{t}=t_{m+1} $.
	
\noindent \textbf{Establishment of Claim 2:} {We prove the claim by contradiction. Suppose that there }exists a time $\bar{t}\in (t_{m-1} ,t_{m+1} )$ such that
\begin{equation}\label{claim2:contr}
|z_{m} (\bar{t})|\, =\zeta _{m}\,  \textrm{and}\, |z_{m} (t)|\, <\zeta _{m},\, \forall \, t\in [t_{m-1} ,\bar{t}).
\end{equation}

\noindent Therefore, by taking into account \eqref{switching:obs}, \eqref{switching:obs:def} and \eqref{claim2:contr}, the map $z_{m} (\cdot )$ satisfies:
\begin{align}
\dot{z}_{m}(t)  =&F(t-\tau,z_m(t),y_{\tau}(t),u_{\tau}(t)) \nonumber\\
& +\phi _{m} (t)P_{m}^{-1} (t)H'(t-\tau,u_{\tau}(t))(y_{\tau}(t)-H(t-\tau,u_{\tau}(t))z_m(t)), \forall \, t\in [t_{m-1} ,\bar{t}]. \label{map:zm:obs}
\end{align}

\noindent Define
\begin{equation} \label{error:switching}
e_{m} (t):=x_{\tau}(t)-z_{m}(t),\forall \, t\in [t_{m-1} ,t_{m+1} ].
\end{equation}
	
\noindent According to Claim 1, $Q_{\, m} (\cdot ):=Q_{m,t_{0} ,\xi _{m} } (\cdot )$ satisfies the CP and \eqref{rhs:difference} with $R:=m$, $t_{0}\ge 0 $  and $\xi :=\xi _{m} $. In addition, due to \eqref{switching:properties}, relations \eqref{map:PR:properties}, \eqref{map:dRbar:int} and \eqref{Lyapunov:Ineq:prop} {are also satisfied with $\bar{t}_{0} :=t_{m-1} $, $R=m$ and $t_0$, $\tau$, $\xi$ as above}. Moreover, due to definition \eqref{xi:i}, the constant $\xi (=\xi _{m} )$ satisfies \eqref{xi} with $R=m$ and $\bar{t}_{0} =t_{m-1} $. Finally, since $|x_{0} |\, \le m=R$, we have $(u,y)\in O(t_0,B_{m} \cap M)$. The previous properties, in conjunction with the result of Proposition \ref{Proposition:state_deter}(ii), assert that estimation {\eqref{error:estimate:2} holds}, namely:
\begin{equation}\label{error:s:estimation:2}
|e_{m} (t)|\, <{\beta (t_{m-1} ,m)}\sqrt{L} \exp \left(-k_{m}(t)+1 \right), \forall \, t\in [t_{m-1} ,\bar{t}].
\end{equation}

\noindent By taking into account {\eqref{int:s:dm}, \eqref{xi:i} and \eqref{error:s:estimation:2}} we get:
\begin{equation} \label{error:s:bound}
|e_{m} (t)|<\xi _{m} , \forall  t\in [t_{m-1} ,\bar{t}].
\end{equation}
	
\noindent Then, from \eqref{state:bound},\eqref{switching:obs:zeta}, \eqref{error:switching}, \eqref{error:s:bound} and by taking into account that $m\ge m_{0} $, {$\beta\in NN$ and $\bar{t}\le t_{m+1} $ we deduce:
\begin{equation} \label{zeta:s:bound}
|z_{m} (\bar{t})|\, \le \, |x(\bar{t})|+|e_{m} (\bar{t})|\, <\beta (\bar{t},m)+\xi _{m} \le \beta (t_{m+1},m)+\xi _{m} =\zeta _{m},
\end{equation}}

\noindent which contradicts \eqref{claim2:contr}. We conclude that for the above $m$, $\bar{t}=t_{m+1} $ and \eqref{claim2:prop2} holds. The same arguments above also assert that \eqref{error:s:estimation:2} holds for every $t\in [t_{m-1} ,t_{m+1} ]$. Notice next that \eqref{int:s:dm} and \eqref{recursive:time} imply that
\begin{equation}
{\exp \left(-\int _{t_{m-1} }^{t}\bar{d}_{m} (s)ds \right)\le \exp(-\kappa_{m}(t)+1)\le\frac{1}{m \beta(t_{m-1},m)\sqrt{L}},\forall t\ge t_m},\label{int:s:estim}
\end{equation}
	
\noindent therefore, by taking into account \eqref{ineq:req:time}, \eqref{error:s:estimation:2} and \eqref{int:s:estim} we have:
\begin{equation} \label{error:s:estim:final}
|e_{m} (t)|\, <\beta (t_{m-1},m)\sqrt{L} \exp(-\kappa_{m}(t)+1)\le \frac{1}{m} , \forall \, t\in [t_{m} ,t_{m+1} ].
\end{equation}
	
\noindent Finally, we show that $\mathop{\lim }\limits_{t\to \infty } |x_{\tau}(t)-Z(t)|\, =0$ or equivalently:
\begin{equation} \label{error:lim:0}
\forall \varepsilon >0 \Rightarrow \exists \, T>0: |x_{\tau}(t)-Z(t)|\, <\varepsilon , \forall t>T.
\end{equation}

\noindent Indeed, let $\varepsilon >0$, $\bar{m}=\bar{m}(\varepsilon )\in {\mathbb N}$ with $\bar{m}\ge m_{0} $ such that $\frac{1}{\bar{m}} <\varepsilon $ and $T=T(\varepsilon ):=t_{\bar{m}} $. Then by \eqref{time:prop} it follows that for every $t>T$ there exists $m\ge \bar{m}$ such that $t_{m} \le t<t_{m+1} $, therefore, since $m\ge m_{0} $ and by recalling {\eqref{error:switching}, \eqref{error:s:estim:final} and that $Z(t):=z_{m} (t)$, $t\in [t_{m} ,t_{m+1} )$}, it follows that $|x_{\tau}(t)-Z(t)|\, =\, |x_{\tau}(t)-z_{m} (t)|\, =\, |e_{m} (t)|\, \le \frac{1}{m} \le \frac{1}{\bar{m}} <\varepsilon $ for every $t>T=T(\varepsilon )$. The latter implies \eqref{error:lim:0} and the proof of Proposition \ref{Proposition:switching} is completed.
\end{IEEEproof}

\section{Proof of Proposition \ref{Theorem}.}\label{section:main_result}

In this section we prove our main result concerning the solvability of the DSODP(DODP) for triangular control systems \eqref{system:triangular}. {The proof of both statements of Proposition \ref{Theorem} is based on the results of Propositions \ref{Proposition:state_deter} and \ref{Proposition:switching}, and partially extends the methodology of the main result in \cite{BdTj13a}. }

\begin{IEEEproof}[Proof of Proposition \ref{Theorem}]
The proof of the first statement, is based on the establishment of Hypothesis \ref{hypothesis:switchin} for system \eqref{system:linearout}. {Notice first that system \eqref{system:triangular} has the form \eqref{system:linearout} with}
\begin{align} 
	F(t,x,y,u):=&\left(\begin{array}{c} f_{1}(t,y,u)+a_{1}(t,y,u)x_{2} \\ f_{2}(t,y,x_{2},u)+a_{2}(t,y,u)x_{3} \\ \vdots \\ f_{n-1}(t,y,x_{2},\ldots,x_{n-1},u)+a_{n-1}(t,y,u)x_{n} \\ f_{n}(t,y,x_{2},\ldots,x_{n},u) \end{array} \right),\nonumber\\
	&(t,x,y,u)\in \RgeO\times\Rat{n}\times\Rat{}\times\Rat{p},\label{rhs:var:triang}
\end{align}
and
\begin{equation} \label{output:triang}
	H:=(\underbrace{1,0,\ldots,0}_{n}).
\end{equation}

\noindent {We show that there exist an integer $\ell\in {\mathbb N}$, a continuous mapping $A:\RgeO\times\Rat{n}\times\Rat{}\times\Rat{p}\to\Rat{n\times n}$ and a constant $L>1$, in such a way that for each $R>0$, both A1 and A2 of Hypothesis \ref{hypothesis:det} hold for system \eqref{system:linearout} with $F(\cdot,\cdot,\cdot,\cdot)$ and $H$ as above}. Let $R>0$ and $\xi>0$. By the $C^{1}$ assumption on the dynamics of {\eqref{system:triangular:st}} and by taking into account \eqref{u:bound} of hypothesis H1, there exists a continuously differentiable function $\sigma_{R}:=\sigma_{R,\xi}\in N$ satisfying:
\begin{align} \label{map:sigma}
	\sigma_{R}(t)\ge\Biggl[&\sum_{i=2}^{n}\sum_{j=2}^{i}\Bigl(\max\lbrace\left|\frac{\partial f_{i}}{\partial x_{j}}(t,x_{1},\ldots,x_{i},u)\right|:|(x_{1},\ldots,x_{i},u)|\le 2\beta(t,R)+\bar{u}(t)+\xi\rbrace\Bigr)^{2}\Biggr]^{\frac{1}{2}},t\ge 0.
\end{align}
	
\noindent Next, consider the set-valued map $[0,\infty)\ni t\to Q_{R}(t):=Q_{R,\xi}(t)\subset {\mathbb R}^{\ell }$, $\ell:=\frac{n(n+1)}{2}$ defined as
\begin{align} \label{map:Q:triang}
	Q_{R}(t):=\{ q=(q_{1,1};q_{2,1},q_{2,2};\ldots;q_{n,1},q_{n,2},&\ldots,q_{n,n})\in {\mathbb R}^{\ell}:|q|\le\sigma_{R}(t)\},
\end{align}
	
\noindent that obviously satisfies the CP. Also, let $Y_{R}(\cdot)$ as given by \eqref{output:set}, with $H(\cdot)$ as given by \eqref{output:triang}, and notice that, due to  \eqref{output:set} and \eqref{output:triang}, it holds:
\begin{equation} \label{output:bound}
	|y|\le\beta(t,R),\: {\rm for}\: {\rm every}\: y\in Y_{R}(t), t\ge 0.
\end{equation}
	
\noindent From \eqref{rhs:var:triang}, \eqref{map:sigma}, \eqref{map:Q:triang} and \eqref{output:bound}, it follows that for every $t\ge 0$, $u\in \Rat{p}$ with $|u|\le\bar{u}(t)$, $y\in Y_{R}(t)$ and $x,z\in {\mathbb R}^{n}$ with $|x|\le\beta(t,R)$ and $|x-z|\le\xi$ we have:
\begin{subequations}
	\begin{align}
		&F(t,x,y,u)-F(t,z,y,u)=A(t,q,y,u)(x-z),\nonumber \\
		&{\rm for}\: {\rm some}\: q\in Q_{R}(t)\: {\rm with}\: q_{i,1}=0,i=1,\ldots,n; \label{rhs:difference:triang}
		\end{align}
		
\noindent with
	\begin{align}
		A(t,q,y,u):=\left(\begin{matrix}
		q_{1,1} & a_{1}(t,y,u) & 0 & \cdots & 0 \\
		q_{2,1} & q_{2,2} & a_{2}(t,y,u) & \ddots & \vdots \\
		\vdots & \vdots & \vdots & \ddots & 0 \\
		q_{n-1,1} & q_{n-1,2} & q_{n-1,3} &  & a_{n-1}(t,y,u) \\
		q_{n,1} & q_{n,2} & q_{n,3} & \cdots & q_{n,n} \end{matrix}\right). \label{map:A:triang}
	\end{align}
\end{subequations}
	
\noindent Hence, A1 is satisfied. In order to establish A2, we prove that there exists a constant $L>1$, such that for every $R>0$, $\xi>0$, $\tau>0$, $t_0\ge 0$,  $\bar{t}_{0}\ge t_{0}+\tau$ and $(u,y)\in O(t_{0},M)$, there exist a nondecreasing continuous function $\kappa_{R}:=\kappa_{R,\xi,\tau,\bar{t}_{0}}\in C^{0}([\bar{t}_{0},\infty);\Rat{})$ satisfying \eqref{map:kappaR} and a sequence of intervals $\mathcal{A}_{\nu}=[\alpha_{\nu},\beta_{\nu}]$, ${\nu\in\mathbb{N}}$, in such a way that  a time-varying symmetric matrix $P_{R}:=P_{R,\xi,t_{0},\bar{t}_{0},\tau,y,u}\in C^{1}([\bar{t}_{0},\infty);\Rat{n\times n})$ and a piecewise continuous function $d_{R}:=d_{R,\xi,t_{0},\bar{t}_{0},\tau,y,u}:[\bar{t}_{0},\infty)\to\Rat{}$ can be found, both {causal and strongly} causal on each $\mathcal{A}_{\nu}$ with respect to $O(t_{0},M)$ and satisfying \eqref{map:PR:properties}-\eqref{Lyapunov:ker:inequality},  with $H$ and $Q_{R}(\cdot)$ as given in \eqref{output:triang} and \eqref{map:Q:triang}, respectively. We proceed as follows. Pick $L>1$, and let $R>0$, $\xi>0$ and $\tau>0$. Consider the $C^{1}$ nondecreasing functions $\sigma_{R,i}(\cdot)$, $i=1,\ldots,n-1$ satisfying
\begin{subequations} \label{maps:sigmaRi:sigmaRbar}
	\begin{align} 
		\sigma_{R,i}(t):=&\sigma_{R,\xi,i}(t)\ge\max\{|a_{i}(\bar{t},y,u)|:\bar{t}\in[0,t],|y|\le\beta(\bar{t},R),|u|\le\bar{u}(t)\},\forall t\ge 0,i=1,\ldots,n-1\label{maps:sigmaRi}
	\end{align}
		
\noindent and define:
	\begin{equation} \label{map:sigmaRbar}
		\bar{\sigma}_{R}(t):=\bar{\sigma}_{R,\xi}(t)=\left(\sigma_{R}^{2}(t)+\sum_{i=1}^{n-1}\sigma_{R,i}^{2}(t)\right)^{\frac{1}{2}},t\ge 0
	\end{equation}
\end{subequations}
	
\noindent with $\sigma_{R}(\cdot)$ as given by \eqref{map:sigma}. It then follows from \eqref{map:sigma}, \eqref{map:Q:triang}, \eqref{rhs:difference:triang}, \eqref{map:A:triang}, \eqref{maps:sigmaRi}, \eqref{map:sigmaRbar} and the fact that $\bar{\sigma}_{R}\in N$, that the map $A(\cdot,\cdot,\cdot,\cdot)$ satisfies
\begin{align} 
	&|A(t-\tau,q,y,u)|\le |A(t-\tau,q,y,u)|_{F}\le \bar{\sigma}_{R}(t),\nonumber\\
	&\forall t\ge\tau,q\in Q_{R}(t-\tau),y\in Y_{R}(t-\tau), |u|\leq \bar{u}(t-\tau)\label{map:A:triang:Frobenius:norm}
\end{align}


\noindent By taking into account \eqref{persistency:condition}, we can assume that without any loss of generality, for each $\nu\in\mathbb{N}$, the difference $T_{\nu}-T_{\nu-1}$ is an integer multiple of the delay $\tau$, namely, it holds $T_{\nu}-T_{\nu-1}={j}_{\nu}\tau$ for certain ${j}_{\nu}\in\mathbb{N}$. Also, we assume that without any loss of generality, it holds $\bar{t}_{0}=t_{0}+\tau(=T_{0}+\tau)$ and select a nondecreasing function $\kappa_{R}\in C^{0}([\bar{t}_{0},\infty);\Rat{})$ as:
\begin{equation} \label{map:kappaR:dfn}
	\kappa_{R}(t)\left\lbrace \begin{array}{ll} \vspace{0.15cm}
	:=-\int_{\bar{t}_{0}}^{\bar{t}_{0}+T_{1}+\tau}\bar{\sigma}_{R}(s)ds-1, & t\in [T_{0}+\tau,T_{1}+\tau] \\ \vspace{0.15cm}
	\in\left[-\int_{\bar{t}_{0}}^{\bar{t}_{0}+T_{1}+\tau}\bar{\sigma}_{R}(s)ds-1,1\right], & t\in\left[T_{1}+\tau,\frac{T_{1}+T_{2}}{2}+\tau\right] \\ \vspace{0.15cm}
	:=\nu-1, & t\in\left[\frac{T_{\nu-1}+T_{\nu}}{2}+\tau,T_{\nu}+\tau\right], \, \,\nu=2,3,\ldots \\
	\in [\nu-1,\nu], & t\in\left[T_{\nu}+\tau,\frac{T_{\nu}+T_{\nu+1}}{2}+\tau\right],\, \, \nu=2,3,\ldots
	\end{array}\right.
\end{equation}
	
\noindent We proceed by showing that for every $(u,y)\in O(t_{0},M)$ there exist mappings $P_{R}(\cdot)$ and $d_{R}(\cdot)$ satisfying the desired causality properties and in such a way that \eqref{map:PR:properties}-\eqref{Lyapunov:ker:inequality} hold, with $\kappa_{R}(\cdot)$ as given in \eqref{map:kappaR:dfn}. Define
\begin{equation} \label{integer:jnu}
	j_{\nu}:=\frac{{T_{\nu}-T_{\nu-1}}}{\tau}\in\mathbb{N}, \nu\in\mathbb{N}.
\end{equation}
	
\noindent  Also, let $(u,y)\in O(t_{0},M)$ and define:
\begin{align}
		m_{\nu}&:=\min\{j\in\{1,\ldots,j_{\nu}\}:\exists t\in (T_{\nu-1}+\tau_{}(j-1),T_{\nu-1}+\tau_{} j):a_{i}(t,y(t),u(t))\ne 0,\forall i=1,\ldots,n-1\}, \label{integer:mnu}
\end{align}
	
\noindent which by virtue of \eqref{persistency:condition} is well defined. In particular, having partitioned each interval $[T_{\nu-1},T_{\nu}]$ into $j_{\nu}$  subintervals of length $\tau$, the $m_{\nu}$-th is the first of the subintervals for which $a_{i}(t,y(t),u(t))\neq 0$ for all $i=1,\ldots, n-1$.

\noindent The existence of the desired mappings $P_{R}(\cdot)$ and $d_{R}(\cdot)$ is based on the establishment of the following claim.

\noindent \textbf{Claim 1.} For $u(\cdot)$, $y(\cdot)$ as above and for each $\nu\in\mathbb{N}$, there exist constants $0<\varepsilon_{2,\nu}<\cdots<\varepsilon_{n,\nu}$ with
\begin{equation}
	\varepsilon_{m,\nu}\le \frac{m}{4n^{2}},m=2,\ldots,n
\end{equation}
	
\noindent and intervals
\begin{equation} \label{intervals:Inu:Jnu}
	I_{\nu}:=[a_{\nu},b_{\nu}]\subset J_{\nu}:=[\bar{a}_{\nu},\bar{b}_{\nu}]\subset(T_{\nu-1}+\tau m_{\nu},T_{\nu-1}+\tau(m_{\nu}+1)), m=2,\ldots,n
\end{equation}
	
\noindent {with $\bar{a}_{\nu}<a_{\nu}<b_{\nu}<\bar{b}_{\nu}$ and} $m_{\nu}$ as given in \eqref{integer:mnu}, such that
\begin{align}
	&a_{i}(t-\tau,y_{\tau}(t),u_{\tau}(t))\ne 0,\forall t\in J_{\nu},i=1,\ldots,n-1 \label{ais:ne0:inInu} \\
	&\bar{a}_{\nu}<a_{\nu}-\varepsilon_{m,\nu}<b_{\nu}+\varepsilon_{m,\nu}<\bar{b}_{\nu},m=2,\ldots,n \label{epsilonm:intervals}
\end{align}

\noindent (In Fig. \ref{fig:intervals} we depict the time instants $T_{\nu}$, $\nu\in\mathbb{N}$ (given by Hypothesis H2), as well as the intervals $I_{\nu}$ and $J_{\nu}$ defined in \eqref{intervals:Inu:Jnu}). In addition, for each $m=2,\ldots,n$, there exist mappings $d_{R,m}\in C^1(\cup_{\nu\in\mathbb{N}}J_{\nu};\Rat{})$ and $P_{R,m}\in C^1(\cup_{\nu\in\mathbb{N}}J_{\nu};\Rat{n\times n})$, {both strongly causal on each $J_{\nu}$, $\nu\in\mathbb{N}$ with respect to $O(t_0,M)$}, such that for each $\nu\in\mathbb{N}$ the following hold:

\begin{figure}
\centering\resizebox{!}{4.5cm}{%
\begin{tikzpicture}
\begin{scope}[line width=0.6pt]
\draw(1,3)--(13.5,3);
\draw (1,2.8)--(1,3.2);\draw (1.7,2.9)--(1.7,3.1);\draw (2.4,2.9)--(2.4,3.1);\draw (3.1,2.9)--(3.1,3.1);
\draw (3.8,2.9)--(3.8,3.1);\draw (4.5,2.9)--(4.5,3.1);\draw (5.2,2.9)--(5.2,3.1);\draw (5.9,2.9)--(5.9,3.1);\draw (6.6,2.9)--(6.6,3.1);
\draw (7.3,2.9)--(7.3,3.1);\draw (10,2.9)--(10,3.1);\draw (10.7,2.9)--(10.7,3.1);\draw (11.4,2.9)--(11.4,3.1);\draw (12.1,2.9)--(12.1,3.1);\draw (12.8,2.9)--(12.8,3.1);

\draw[<->] (1,2.85)--(1.7,2.85);
\end{scope}

\draw[dotted] (1.81,1)--(1.81,3);\draw[dotted] (2.31,1)--(2.31,3);
\draw[dashed] (3.8,1)--(3.8,3);
\draw[dotted] (5.3,1)--(5.3,3);\draw[dotted] (5.8,1)--(5.8,3);
\draw[dashed] (1,1)--(1,3);
\draw[dashed] (7.3,1)--(7.3,3);
\draw[dashed] (10,1)--(10,3);
\draw[dashed] (12.8,1)--(12.8,3);
\draw[dotted] (11.45,1)--(11.45,3);
\draw[dotted] (12.05,1)--(12.05,3);

\path [decoration={zigzag,segment length=2.3,amplitude=1.3,
  post=lineto},font=\scriptsize,
  line join=round] (1.81,3) edge[decorate] (2.31,3);

\path [decoration={zigzag,segment length=2.3,amplitude=1.3,
  post=lineto},font=\scriptsize,
  line join=round] (2.9,3) edge[decorate] (3.7,3);  
  
\path [decoration={zigzag,segment length=2.3,amplitude=1.3,
  post=lineto},font=\scriptsize,
  line join=round] (5.3,3) edge[decorate] (6.8,3);
  
\path [decoration={zigzag,segment length=2.3,amplitude=1.3,
  post=lineto},font=\scriptsize,
  line join=round] (11.45,3) edge[decorate] (12.8,3); 
  

\path [decoration={zigzag,segment length=2.3,amplitude=1.3,
  post=lineto},font=\scriptsize,
  line join=round] (5,3.5) edge[decorate] (6.44,3.5);  
  
\node at (9.41,3.5) {:  $a_{i}(\cdot,y(\cdot),u(\cdot))\neq 0,\,\,i=1,\ldots,n-1$};
  
\begin{scope}[line width=0.6pt]
\draw(1,1)--(12.8,1);\draw[dotted](12.8,1)--(13.5,1);
\draw (1,0.8)--(1,1.2);
\draw (3.8,0.8)--(3.8,1.2);\draw (7.3,0.8)--(7.3,1.2);\draw (10,0.8)--(10,1.2);\draw (12.8,0.8)--(12.8,1.2);
\end{scope}

\begin{scope}[line width=0.6pt]
\draw (1.81,0.9)--(1.81,1.1);\draw (2.31,0.9)--(2.31,1.1); 
\draw[<->](1.81,1.2)--(2.31,1.2);

\draw (5.3,0.9)--(5.3,1.1);\draw (5.8,0.9)--(5.8,1.1);
\draw[<->](5.3,1.2)--(5.8,1.2);

\draw (11.45,0.9)--(11.45,1.1);\draw (12.05,0.9)--(12.05,1.1);
\draw[<->](11.45,1.2)--(12.05,1.2);
\end{scope}

\node[text height=0.2cm] at (1,0.5) {$T_0=t_0$};
\node[text height=0.2cm] at (3.8,0.5) {$T_1$};
\node[text height=0.2cm] at (7.34,0.5) {$T_2$};
\node[text height=0.2cm] at (10,0.5) {$T_{\nu-1}$};
\node[text height=0.2cm] at (12.8,0.5) {$T_{\nu}$};
\node[text height=0.3cm] at (8.5,0.9) {$\cdots$};
 \node[text height=0.3cm] at (2.1,0.8) {$J_{1}$};
 \node[text height=0.3cm] at (5.6,0.8) {$J_{2}$};
 \node[text height=0.3cm] at (11.75,0.8) {$J_{\nu}$};
 \node[text height=0.3cm] at (1.35,2.75) {$\tau$};
 \node[text height=0.3cm] at (13.3,3.35) {$t$};
 \node[text height=0.3cm] at (1,3.55) {$t_0$};
 
 \draw[semithick] (11.2,1.6) rectangle (12.3,0.5);
\path [draw,semithick,fill=gray!20] (2.3,0) rectangle (12.6,-1.5);

\begin{scope}[line width=0.6pt]
\draw(2.6,-0.7)--(12.4,-0.7);
\draw (3,-0.5)--(3,-0.9);\draw (5.5,-0.5)--(5.5,-0.9);\draw (9.5,-0.5)--(9.5,-0.9);\draw (12,-0.5)--(12,-0.9);
\end{scope}

\draw(5,-0.55)--(5,-0.85);\node[text height=0.45cm] at (5,-0.85) {$\varepsilon_1$};
\draw(4.5,-0.55)--(4.5,-0.85);\node[text height=0.45cm] at (4.5,-0.85) {$\varepsilon_2$};
\node[text height=0.75cm] at (4,-0.7) {$\cdots$};
\draw(3.5,-0.55)--(3.5,-0.85);\node[text height=0.45cm] at (3.5,-0.85) {$\varepsilon_n$};

\draw(10,-0.55)--(10,-0.85);\node[text height=0.45cm] at (10,-0.85) {$\varepsilon_1$};
\draw(10.5,-0.55)--(10.5,-0.85);\node[text height=0.45cm] at (10.5,-0.85) {$\varepsilon_2$};
\node[text height=0.75cm] at (11,-0.7) {$\cdots$};
\draw(11.5,-0.55)--(11.5,-0.85);\node[text height=0.45cm] at (11.5,-0.85) {$\varepsilon_n$};

\node[text height=0.5cm] at (3,-0.15) {$\bar{a}_{\nu}$};
\node[text height=0.5cm] at (5.5,-0.15) {${a}_{\nu}$};
\node[text height=0.5cm] at (9.5,-0.15) {$b_{\nu}$};
\node[text height=0.5cm] at (12,-0.15) {$\bar{b}_{\nu}$};
\draw[<->](5.5,-0.9)--(9.5,-0.9);
\node[text height=0.5cm] at (7.5,-1.1) {$I_{\nu}$};

\draw[semithick] (11.8,0.5)--(7.5,0);
\end{tikzpicture} }
\caption{Representation of the time instants $T_{\nu}$ and the intervals $I_{\nu}$ and $J_{\nu}$.} \label{fig:intervals}
\end{figure}
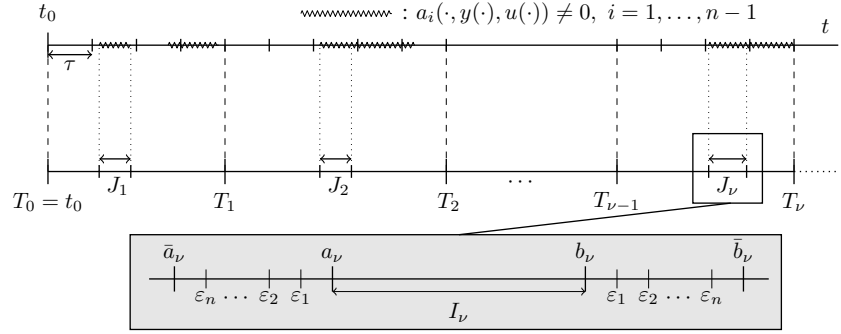

\begin{equation}
\left\{\begin{array}{l}
	P_{R,2}(t)=\left(\begin{matrix}
	p_{R,2,1}(t) & p_{R,2}(t) \\ p_{R,2}(t) & L
	\end{matrix}\right), \,\, \forall t\in J_{\nu}, \\[1.5em] 
	P_{R,m}(t)=\left(\begin{matrix}
	p_{R,m,1}(t) & \begin{matrix} p_{R,m}(t) & 0 & \cdots & 0 \end{matrix} \\ 
	\begin{matrix} p_{R,m}(t) \\ 0 \\ \vdots \\ 0 \end{matrix} & \boxed{P_{R,m-1}(t)}
	\end{matrix}\right), \,\,\,\forall t\in J_{\nu},m=3,\ldots,n,\\[0.5em]
	P_{R,m}(t)>I_{{m\times m}}, \,\,\, \forall t\in J_{\nu},m=2,\ldots,n, \\[0.5em]
	p_{R,m,1}(t)=L,p_{R,m}(t)=0, \,\,\,\forall t\in J_{\nu}\setminus [a_{\nu}-\varepsilon_{m,\nu},b_{\nu}+\varepsilon_{m,\nu}],m=2,\ldots,n
\end{array}\right. \label{map:PRm:properties}
\end{equation}
\begin{equation}
\left\{\begin{array}{l}
	d_{R,m}(t)=\frac{\int_{\bar{t}_{0}}^{T_{2}+\tau}\bar{\sigma}_{R}(s)ds+2}{b_{1}-a_{1}}+n-m,\,\,\, \forall t\in I_{1} \;{\rm if}\;\nu=1,\vspace{0.23em}\\ [0.6em]
	d_{R,m}(t)=\frac{\int_{T_{\nu}+\tau}^{T_{\nu+1}+\tau}\bar{\sigma}_{R}(s)ds+2}{b_{\nu}-a_{\nu}}+n-m, \,\,\,\forall t\in I_{\nu},m=2,\ldots,n, \;{\rm if}\;\nu=2,3,\ldots \\[0.6em]
	\int_{S}d_{R,m}(s)ds\ge -\frac{m}{n}, \,\,\, \forall S\subset [a_{\nu}-\varepsilon_{m,\nu},b_{\nu}+\varepsilon_{m,\nu}]\setminus I_{\nu},\,\,\,m=2,\ldots,n \\[0.6em]
	d_{R,m}(t)\le -\bar{\sigma}_{R}(t), \,\,\, \forall t\in J_{\nu}\setminus [a_{\nu}-\varepsilon_{m,\nu},b_{\nu}+\varepsilon_{m,\nu}],\,\,\,m=2,\ldots,n 
\end{array}\right.\label{map:dRm:properties}
\end{equation}
\begin{align}
	e'P_{R,m}(t)A_{m}(t-\tau,q,y_{\tau}(t),u_{\tau}&(t))e+\tfrac{1}{2}e'\dot{P}_{R,m}(t)e\le-d_{R,m}(t)e'P_{R,m}(t)e, \nonumber \\ \forall t\in J_{\nu}, e\in\ker H_{m},q\in &Q_{R}(t-\tau) \label{Lyapunov:inequality:m}
\end{align}
	
\noindent with
\begin{equation} \label{map:Hm}
	H_{m}:=(\underbrace{1,0,\ldots,0}_{m}),e:=(e_{n-m+1},\ldots,e_{n})'\in {\mathbb R}^{m}
\end{equation}
\begin{align}
	A_{m}(t,q,y,u):=
	&\left(\begin{matrix} q_{n-m+1,n-m+1} & \begin{matrix} a_{n-m+1}(t,y,u) & 0 & \cdots & 0 \end{matrix} \\ \begin{matrix} q_{n-m+2,n-m+1} \\ \vdots \\  q_{n,n-m+1} \end{matrix} & \boxed{A_{m-1}(t,q,y,u)} \end{matrix}\right);\nonumber\\
	&A_{1}(t,q,y,u):=q_{n,n} \label{map:Am}
\end{align}

\begin{IEEEproof}[Proof of Claim 1] 
\textbf{{Case} $m:=2$:} {Consider the constants $L>1$, $R>0$, $\xi>0$, $\tau>0$ as above and define}
\begin{subequations} \label{maps:triangular:2}
	\begin{equation}
	H_{2}:=(1,0),e:=(e_{n-1},e_{n})'\in {\mathbb R}^{2} \label{map:H2}
	\end{equation}
	\begin{equation}
	A_{2}(t,q,y,u):=\left(\begin{matrix} q_{n-1,n-1} & a_{n-1}(t,y,u) \\ q_{n,n-1} & q_{n,n} \end{matrix}\right). \label{map:A2}
	\end{equation}
\end{subequations}

We establish existence of a time-varying symmetric matrix $P_{R,2}:=P_{R,\xi,\bar{t}_{0},\tau,y,u,2}\in C^{1} (\cup_{\nu\in\mathbb{N}}J_{\nu};$ ${\mathbb R}^{2\times 2})$ and a mapping $d_{R,2}:=d_{R,\xi,\bar{t}_{0},\tau,y,u,2}\in C^{1}(\cup_{\nu\in\mathbb{N}}J_{\nu};{\mathbb R})$ in such a way that {for each $\nu\in \mathbb{N}$}
\begin{subequations}
	\begin{equation} \label{map:PR2:properties}
	\left.\begin{array}{ll} \vspace{0.15cm}
	\bullet\; P_{R,2}(t)=\left(\begin{matrix}
	p_{R,2,1}(t) & p_{R,2}(t) \\
	p_{R,2}(t) & L
	\end{matrix}\right), \,\,\,\,\,\,  \,\,\,\,\,\,\forall t\in J_{\nu}, \\ \vspace{0.15cm}
	\bullet\; P_{R,2}(t)>I_{{2\times 2}},  \,\,\,\,\,\, \,\,\,\,\,\, \forall t\in J_{\nu}, \\ \vspace{0.15cm}
	\bullet\; p_{R,2,1}(t)=L,p_{R,2}(t)=0, \forall t\in J_{\nu}\setminus [a_{\nu}-\varepsilon_{2,\nu},b_{\nu}+\varepsilon_{2,\nu}],
	\end{array}\right.
	\end{equation}
	
	\begin{equation} \label{map:dR2:properties}
	\left.\begin{array}{ll} \vspace{0.15cm}
	\bullet\; d_{R,2}(t)=\frac{\int_{t_{0}}^{T_{2}}\bar{\sigma}_{R}(s)ds+2}{b_{1}-a_{1}}+n-2,  \,\,\,\,\,\, \forall t\in I_{1},\;\;\nu=1  \\ \vspace{0.15cm}
	\bullet\; d_{R,2}(t)=\frac{\int_{T_{\nu}}^{T_{\nu+1}}\bar{\sigma}_{R}(s)ds+2}{b_{\nu}-a_{\nu}}+n-2,  \,\,\,\,\,\, \forall t\in I_{\nu},\;\;\nu\ge 2 \\ \vspace{0.15cm}
	\bullet\; \int_{S}d_{R,2}(s)ds\ge -\frac{2}{n},  \,\,\,\,\,\, \forall S\subset [a_{\nu}-\varepsilon_{2,\nu},b_{\nu}+\varepsilon_{2,\nu}]\setminus I_{\nu}, \\  \vspace{0.15cm}
	\bullet\; d_{R,2}(t)\le -\bar{\sigma}_{R}(t),  \,\,\,\,\,\, \forall t\in[T_{\nu-1},T_{\nu}]\setminus J_{\nu},
	\end{array}\right.
	\end{equation}
	\begin{align}
	e'P_{R,2}(t)A_{2}(t-\tau,q,y_{\tau}(t),u_{\tau}(t))e+&\tfrac{1}{2}e'\dot{P}_{R,2}(t)e\le-d_{R,2}(t)e'P_{R,2}(t)e, \nonumber \\
	\forall t\in J_{\nu},e\in\ker H_{2},q\in Q_{R}&(t-\tau) \label{Lyapunov:inequality:2}
	\end{align}
\end{subequations}

\noindent with $H_{2}$, $A_{2}(\cdot,\cdot,\cdot,\cdot)$ and $Q_{R}(\cdot)$ as given in \eqref{map:H2}, \eqref{map:A2} and \eqref{map:Q:triang}, respectively. By taking into account \eqref{maps:triangular:2} and \eqref{map:PR2:properties}, the desired \eqref{Lyapunov:inequality:2} is written:
\begin{align}
(0,e)&\left(\begin{matrix} p_{R,2,1}(t) & p_{R,2}(t) \\ p_{R,2}(t) & L \end{matrix}\right)\left(\begin{matrix} q_{n-1,n-1} & a_{n-1}(t-\tau,y_{\tau}(t),u_{\tau}(t)) \\ q_{n,n-1} & q_{n,n} \end{matrix}\right)\left(\begin{matrix} 0 \\ e \end{matrix}\right) \nonumber \\
&+(0,e)\dot{\overbrace{\left(\begin{matrix} p_{R,2,1}(t) & p_{R,2}(t) \\ p_{R,2}(t) & L \end{matrix}\right)}}\left(\begin{matrix} 0 \\ e \end{matrix}\right)
\le-d_{R,2}(t)(0,e)\left(\begin{matrix} p_{R,2,1}(t) & p_{R,2}(t) \\ p_{R,2}(t) & L \end{matrix}\right)\left(\begin{matrix} 0 \\ e \end{matrix}\right), \nonumber \\
&\forall t\in J_{\nu}, e\in\Rat{},q\in Q_{R}(t-\tau),\iff \nonumber \\
p_{R,2}(t)&a_{n-1}(t-\tau,y_{\tau}(t),u_{\tau}(t))+Lq_{n,n}\le-Ld_{R,2}(t), \nonumber \\
&\forall t\in J_{\nu}, q\in Q_{R}(t-\tau) \label{Lin2:follows:from1}
\end{align}

\noindent By invoking \eqref{map:Q:triang}, {\eqref{map:A:triang:Frobenius:norm}} and the equivalence between \eqref{Lyapunov:inequality:2} and \eqref{Lin2:follows:from1}, it follows that, in order to prove \eqref{Lyapunov:inequality:2}, it suffices to determine $p_{R,1},p_{R}\in C^{1}(\cup_{\nu\in\mathbb{N}}J_{\nu};{\mathbb R})$ and $d_{R,2} \in C^{1}(\cup_{\nu\in\mathbb{N}}J_{\nu};{\mathbb R})$ in such a way that for each $\nu\in\mathbb{N}$, \eqref{map:PR2:properties} and \eqref{map:dR2:properties} are fulfilled, and further:
\begin{align} 
p_{R,2}(t)a_{n-1}(t-\tau,y_{\tau}(t),&u_{\tau}(t))+L\bar{\sigma}_{R}(t)\le-Ld_{R,2}(t),\nonumber\\
&\forall t\in J_{\nu},\nu\in\mathbb{N}\label{Lin2:follows:from2}
\end{align}

\noindent \textbf{Construction of ${p_{R}(\cdot)}$ and ${d_{R,2}(\cdot)}$:} {For each $\nu\in\mathbb{N}$, let}
\begin{subequations} \label{constants:M2:tau2}
	\begin{align}
	M_{2,\nu}:=&\max\{\bar{\sigma}_{R}(t):t\in J_{\nu}\} \label{constant:M2} \\
	\varepsilon_{2,\nu}:=&\min\left\{\frac{1}{2nM_{2}},\frac{\bar{b}_{\nu}-b_{\nu}}{2},\frac{a_{\nu}-\bar{a}_{\nu}}{2},\frac{1}{2n^{2}}\right\} \label{constant:tau2}
	\end{align}
\end{subequations}

\noindent and define $d_{R,2}\in C^{1}(\cup_{\nu\in\mathbb{N}}J_{\nu};{\mathbb R})$ and $p_{R}\in C^{1}(\cup_{\nu\in\mathbb{N}}J_{\nu};{\mathbb R})$ as follows:
\begin{align}
d_{R,2}(t)&\left\lbrace\begin{array}{ll} \vspace{0.15cm}
:=-\bar{\sigma}_{R}(t),  & t\in[\bar{a}_{\nu},a_{\nu}-\varepsilon_{2,\nu}],\nu\in\mathbb{N} \\ \vspace{0.15cm}
\in\left[-\bar{\sigma}_{R}(t),\frac{\int_{t_{0}}^{T_{2}}\bar{\sigma}_{R}(s)ds+2}{b_{1}-a_{1}}+n-2\right], & t\in [a_{1}-\varepsilon_{2,1},a_{1}] \\ \vspace{0.15cm}
:=\frac{\int_{t_{0}}^{T_{2}}\bar{\sigma}_{R}(s)ds+2}{b_{1}-a_{1}}+n-2, & t\in[a_{1},b_{1}] \\ \vspace{0.15cm}
\in\left[-\bar{\sigma}_{R}(t),\frac{\int_{t_{0}}^{T_{2}}\bar{\sigma}_{R}(s)ds+2}{b_{1}-a_{1}}+n-2\right], & t\in [b_{1},b_{1}+\varepsilon_{2,1}] \\ \vspace{0.15cm}
\in\left[-\bar{\sigma}_{R}(t),\frac{\int_{T_{\nu}}^{T_{\nu+1}}\bar{\sigma}_{R}(s)ds+2}{b_{\nu}-a_{\nu}}+n-2\right], & t\in [a_{\nu}-\varepsilon_{2,\nu},a_{\nu}],\nu=2,3,\ldots \\ \vspace{0.15cm}
:=\frac{\int_{T_{\nu}}^{T_{\nu+1}}\bar{\sigma}_{R}(s)ds+2}{b_{\nu}-a_{\nu}}+n-2, & t\in[a_{\nu},b_{\nu}],\nu=2,3,\ldots \\ \vspace{0.15cm}
\in\left[-\bar{\sigma}_{R}(t),\frac{\int_{T_{\nu}}^{T_{\nu+1}}\bar{\sigma}_{R}(s)ds+2}{b_{\nu}-a_{\nu}}+n-2\right], & t\in [b_{\nu},b_{\nu}+\varepsilon_{2,\nu}],\nu=2,3,\ldots \\ \vspace{0.15cm}
:=-\bar{\sigma}_{R}(t),  & t\in[b_{\nu}+\varepsilon_{2,\nu},\bar{b}_{\nu}],\nu\in\mathbb{N}
\end{array}\right. \label{map:dR2:c2} \\
p_{R,2}(t)&:=L\frac{-d_{R,2}(t)-\bar{\sigma}_{R}(t)}{a_{n-1}(t-\tau,y_{\tau}(t),u_{\tau}(t))},t\in J_{\nu},\nu\in\mathbb{N}. \label{map:pR:keq2:dfn:c2}
\end{align}

\noindent  It follows, by taking into account the definition of $d_{R,2}(\cdot)$ and $p_{R,2}(\cdot)$ above, that \eqref{Lin2:follows:from2} holds and consequently the desired \eqref{Lyapunov:inequality:2} is satisfied for all $\nu \in \mathbb{N}$. {By invoking \eqref{integer:jnu}, \eqref{integer:mnu}, \eqref{intervals:Inu:Jnu}-\eqref{epsilonm:intervals}, \eqref{map:dR2:c2}, \eqref{map:pR:keq2:dfn:c2} and the fact that $\sigma_{R}(\cdot)$ is independent of $(u,y)$, it follows that the $p_{R,2}(\cdot)$ in \eqref{map:dR2:c2} above is strongly causal on each $J_{\nu}$, $\nu\in\mathbb{N}$ with respect to $O(t_0,M)$.}

Finally, we show that $\int_{S}d_{R,2}(s)ds\ge-\tfrac{2}{n}$ for any $S\subset [a_\nu-\varepsilon_{2,\nu},b_{\nu}+\varepsilon_{2,\nu}]\setminus I_{\nu}$. Indeed, by invoking  \eqref{map:sigmaRbar}, \eqref{constant:M2}, \eqref{constant:tau2} and \eqref{map:dR2:c2} we get
\begin{align}
\int_{S}d_{R,2}(s)ds&\ge
-\int_{S}\bar{\sigma}_R(s)ds\nonumber\\
&\ge -\int_{[a_{\nu}-\varepsilon_{m,\nu},b_{\nu}+\varepsilon_{m,\nu}]\setminus I_{\nu}}\bar{\sigma}_{R}(s)ds\nonumber\\
&\ge-2M_{2,\nu}\varepsilon_{2,\nu}\geq-\frac{2}{n}.\label{int:dR2:1}
\end{align}

\noindent {From \eqref{map:dR2:c2} and \eqref{int:dR2:1}} we deduce that all properties of \eqref{map:dR2:properties} are satisfied. Hence, \eqref{map:dRm:properties} is fulfilled for $m=2$.

\noindent \textbf{Construction of ${p_{R,1}(\cdot)}$:} {We now determine a function $p_{R,1}\in C^{1}([t_0,\infty);\mathbb{R})$ such that \eqref{map:PRm:properties} holds. Notice first that, from \eqref{map:dR2:c2} and \eqref{map:pR:keq2:dfn:c2} we have $p_{R,2}(t)=0$ for all $t\in J_{\nu}\setminus[a_{\nu}-\varepsilon_{2,\nu},b_{\nu}+\varepsilon_{2,\nu}]$.} Since $L>1$, we may define
\begin{equation}\label{map:pR1}
p_{R,1}(t):=\frac{p_{R,2}^2(t)}{L-1}+L
\end{equation}

\noindent {which implies that $p_{R,1}(t)=L$ for all $t\in J_{\nu}\setminus[a_{\nu}-\varepsilon_{2,\nu},b_{\nu}+\varepsilon_{2,\nu}]$, $\nu\in\mathbb{N}$. Then, by taking into account \eqref{map:pR:keq2:dfn:c2}, \eqref{map:pR1} and the assumption that $L>1$ we derive that}
\begin{align}
\det(P_{R,2}(t)-I_{2\times 2})&=\det\left(\begin{array}{ll}p_{R,1}(t)-1 &p_{R,2}(t)\\p_{R,2}(t) & L-1\end{array}\right)\nonumber\\
&=(L-1)^2>0, \,\forall t\in J_{\nu}.\label{det:pR2}
\end{align}

\noindent {which in turn implies that $P_{R,2}(t)>I_{2\times 2}$ for all $t\in J_{\nu}$.} {It follows from \eqref{map:pR1} and the causality properties of  $p_{R,2}(\cdot)$ as defined in \eqref{map:pR:keq2:dfn:c2}, that  $p_{R,1}(\cdot)$ is also strongly causal on each $J_{\nu}$, $\nu\in\mathbb{N}$  with respect to $O(t_0,M)$.} Finally, from \eqref{map:pR:keq2:dfn:c2}, \eqref{map:pR1} and \eqref{det:pR2} we conclude that all properties of \eqref{map:PR2:properties} are fulfilled as well. Thus, \eqref{map:PRm:properties} holds with $m=2$.  This completes the proof of Claim 1 for $m=2$.

\noindent \textbf{Proof of Claim 1; general step of induction procedure:} Assume now that Claim 1 is fulfilled for certain integer $m$ with $2\le m<n$. We prove that Claim 1 also holds for $m:=m+1$. Consider the pair $(H,A)$ as given in \eqref{maps:HA:aux} with $H(t,u):=H_{m}$, $A(t,q,y,u):=A_{m}(t,q,y,u)$, $\ell =\frac{n(n+1)}{2}$, $m:=m$, $n:=n$ and $s=k:=1$, where $H_{m}$ and $A_{m}$ are defined by \eqref{map:Hm} and \eqref{map:Am},  respectively. Also, consider the set-valued map  $Q_{R}(\cdot)$ as given in \eqref{map:Q:triang} and the {strongly causal on each $J_{\nu}$, $\nu\in\mathbb{N}$} with respect to $O(t_0,M)$, mappings $d_{R,m}(\cdot)$ and $P_{R,m}(\cdot)$ as defined by  \eqref{map:PRm:properties}, \eqref{map:dRm:properties} in the induction hypothesis, and satisfying \eqref{Lyapunov:inequality:m} with $m:=m$. Finally, consider the function $\bar{d}_{R,m}:=\bar{d}_{R,\xi,t_{0},\tau,y,u,m}$ defined as:
\begin{equation} \label{map:dRbm}
\bar{d}_{R,m}(t):=d_{R,m}(t)-1,\,t\in J_{\nu},\,\nu\in\mathbb{N},
\end{equation}

\noindent \noindent which satisfies $\bar{d}_{R,m}(t)< d_{R,m}(t)$ for all $t\in\cup_{\nu\in\mathbb{N}}J_{\nu}$. It follows that all requirements of Lemma \ref{lemma} are fulfilled and therefore, there exists {a strongly causal on each $J_{\nu}$, $\nu\in\mathbb{N}$} with respect to $O(t_0,M)$ function $\phi_{R,m}:=\phi_{R,\xi,t_{0},\tau,m}\in C^{1}(\cup_{\nu\in\mathbb{N}}J_{\nu};{\mathbb R}_{> 0})$ such that
\begin{align}
e'P_{R,m}(t)A_{m}&(t-\tau,q,y_{\tau}(t),u_{\tau}(t))e+\tfrac{1}{2}e'\dot{P}_{R,m}(t)e\le\phi_{R,m}(t)|H_{m}e|^{2}-\bar{d}_{R,m}(t)e'P_{R,m}(t)e, \nonumber \\
\forall t\in J_{\nu},&\nu\in\mathbb{N},q\in Q_{R}(t-\tau),e\in {\mathbb R}^{m} \label{Lyapunov:inequality:aux:m}
\end{align}

\noindent Furthermore, due to \eqref{map:dRm:properties} and \eqref{map:dRbm}, the map $\bar{d}_{R,m}(\cdot)$ satisfies:
\begin{align} \label{map:dRbm:properties}
\begin{array}{ll} \vspace{0.15cm}
\bullet\; \bar{d}_{R,m}(t)=\frac{\int_{\bar{t}_{0}}^{T_{2}+\tau}\bar{\sigma}_{R}(s)ds+2}{b_{1}-a_{1}}+n-(m+1), \,\,\,\,\,\,\, \forall t\in I_{1} \\ \vspace{0.15cm}
\bullet\; \bar{d}_{R,m}(t)=\frac{\int_{T_{\nu}+\tau}^{T_{\nu+1}+\tau}\bar{\sigma}_{R}(s)ds+2}{b_{\nu}-a_{\nu}}+n-(m+1), \,\,\,\,\,\,\,  \forall t\in I_{\nu},\nu=2,3,\ldots \\ \vspace{0.15cm}
\bullet\; \int_{S}\bar{d}_{R,m}(s)ds\ge-\frac{m}{n}-\frac{1}{2n}, \,\,\,\,\,\,\, \forall S\subset [a_{\nu}-\varepsilon_{m,\nu},b_{\nu}+\varepsilon_{m,\nu}]\setminus I_{\nu},\,\nu\in\mathbb{N} \\ \vspace{0.15cm}
\bullet\; \bar{d}_{R,m}(t)\le -\bar{\sigma}_{R}(t), \,\,\,\,\,\,\, \forall t\in  J_{\nu}\setminus [a_{\nu}-\varepsilon_{m,\nu},b_{\nu}+\varepsilon_{m,\nu}], \, \nu\in\mathbb{N}
\end{array}
\end{align}

In the sequel, we exploit \eqref{Lyapunov:inequality:aux:m} and \eqref{map:dRbm:properties}, in order to establish that Claim 1 is fulfilled for $m=m+1$. Specifically, for the same $L$, $R$, $\xi$, $t_{0}=\bar{t}_{0}$, $\tau$, $u(\cdot)$ and $y(\cdot)$ as above, we show that there exist a time-varying symmetric matrix $P_{R,m+1}\in C^{1}(\cup_{\nu\in\mathbb{N}}J_{\nu};{\mathbb R}^{(m+1)\times (m+1)})$ and a map $d_{R,m+1}\in C^{1}(\cup_{\nu\in\mathbb{N}}J_{\nu};{\mathbb R})$, such that both \eqref{map:PRm:properties} and \eqref{map:dRm:properties} are fulfilled with $m=m+1$ and further:
\begin{align}
&e'P_{R,m+1}(t)A_{m+1}(t-\tau,q,y_{\tau}(t),u_{\tau}(t))e+\tfrac{1}{2}e'\dot{P}_{R,m+1}(t)e\le -d_{R,m+1}(t)e'P_{R,m+1}(t)e, \nonumber \\
&\forall t\in J_{\nu},\nu\in\mathbb{N},e\in\ker H_{m+1},q\in Q_{R}(t-\tau) \label{Lyapunov:inequality:mpl1}
\end{align}

\noindent with
\begin{subequations} \label{maps:triangular:mpl1}
	\begin{equation}
	H_{m+1}:=(\underbrace{1,0,\ldots,0}_{m+1}),e:=(e_{n-m};\hat{e}')'\in {\mathbb R}\times {\mathbb R}^{m},\hat{e}:=(e_{n-m+1},\ldots,e_{n})'\in {\mathbb R}^{m} \label{map:Hmpl1}
	\end{equation}
	\begin{equation}
	A_{m+1}(t,q,y,u):=\left(\begin{matrix} q_{n-m,n-m} & \begin{matrix} a_{n-m}(t,y,u) & 0 & \cdots & 0 \end{matrix} \\ \begin{matrix} q_{n-m+1,n-m} \\ \vdots \\  q_{n,n-m} \end{matrix} & \boxed{A_{m}(t,q,y,u)} \end{matrix}\right) \label{map:Ampl1}
	\end{equation}
	\begin{equation}
	P_{R,m+1}(t):=\left(\begin{matrix} p_{R,m+1,1}(t) & \begin{matrix} p_{R,m+1}(t) & 0 & \cdots & 0 \end{matrix} \\ \begin{matrix} p_{R,m+1}(t) \\ 0 \\ \vdots \\ 0 \end{matrix} & \boxed{P_{R,m}(t)} \end{matrix}\right) \label{map:PRmpl1}
	\end{equation}
\end{subequations}

\noindent and where $Q_{R}(\cdot)$ is given in \eqref{map:Q:triang}. Notice that, according to \eqref{map:Hmpl1} and \eqref{map:PRmpl1}, we have $e'P_{R,m+1}(t)e=\hat{e}'P_{R,m}(t)\hat{e}$ for every $e=(0,\hat{e}')'=(0,e_{n-m+1},\ldots,e_{n})'\in \ker H_{m+1}$.  Thus, by taking into account \eqref{map:Hm} and \eqref{maps:triangular:mpl1}, the desired \eqref{Lyapunov:inequality:mpl1} is equivalently written:
\begin{align}
(0;\hat{e}')&\left(\begin{matrix} p_{R,m+1,1}(t) & \begin{matrix} p_{R,m+1}(t) & 0 & \cdots & 0 \end{matrix} \\ \begin{matrix} p_{R,m+1}(t) \\ 0 \\ \vdots \\ 0 \end{matrix} & \boxed{P_{R,m}(t)} \end{matrix}\right)\nonumber\\
&\hspace{8em}\times \left(\begin{matrix} q_{n-m,n-m} & \begin{matrix} a_{n-m}(t-\tau,y_{\tau}(t),u_{\tau}(t)) & 0 & \cdots & 0 \end{matrix} \\ \begin{matrix} q_{n-m+1,n-m} \\ \vdots \\ q_{n,n-m} \end{matrix} & \boxed{A_{m}(t-\tau,q,y,u)} \end{matrix}\right)\left(\begin{matrix} 0 \\ \hat{e} \end{matrix}\right) \nonumber \\
+&\frac{1}{2}(0;\hat{e}')\dot{\overbrace{\left(\begin{matrix} p_{R,m+1,1}(t) & \begin{matrix} p_{R,m+1}(t) & 0 & \cdots & 0 \end{matrix} \\ \begin{matrix} p_{R,m+1}(t) \\ 0 \\ \vdots \\ 0 \end{matrix} & \boxed{P_{R,m}(t)} \end{matrix}\right)}}\left(\begin{matrix} 0 \\ \hat{e} \end{matrix}\right) \nonumber \\
\le& -d_{R,m+1}(t)(0;\hat{e}')\left(\begin{matrix} p_{R,m+1,1}(t) & \begin{matrix} p_{R,m+1}(t) & 0 & \cdots & 0 \end{matrix} \\ \begin{matrix} p_{R,m+1}(t) \\ 0 \\ \vdots \\ 0 \end{matrix} & \boxed{P_{R,m}(t)} \end{matrix}\right)\left(\begin{matrix} 0 \\ \hat{e} \end{matrix}\right), \nonumber \\
&\forall t\in J_{\nu},\nu\in\mathbb{N},\hat e\in\Rat{m},q\in Q_{R}(t-\tau) \iff \nonumber
\end{align}

\begin{align}
e_{n-m+1}^{2}p_{R,m+1}(t)a_{n-m}(t-\tau,y_{\tau}(t),u_{\tau}(t))&+\hat{e}'P_{R,m}(t)A_{m}(t-\tau,q,y_{\tau}(t),u_{\tau}(t))\hat{e}\nonumber\\
&+\tfrac{1}{2}\hat{e}'P_{R,m}(t)\hat{e}\le -d_{R,m+1}(t)\hat{e}'P_{R,m}(t)\hat{e}, \nonumber \\
\forall t\in J_{\nu},\nu\in\mathbb{N},\hat{e}\in\Rat{m},& q\in Q_{R}(t-\tau), \label{Linmpl1:follows:from1}
\end{align}

\noindent  where $A_{m}(\cdot,\cdot,\cdot,\cdot)$ and $P_{R,m}(\cdot)$ are given by \eqref{map:Am} and \eqref{map:PRm:properties}, respectively. From \eqref{Lyapunov:inequality:aux:m}, it suffices, instead of \eqref{Linmpl1:follows:from1}, to show that
\begin{align}
&e_{n-m+1}^{2}(p_{R,m+1}(t)a_{n-m}(t-\tau,y_{\tau}(t),u_{\tau}(t))+\phi_{R,m}(t))\le (\bar{d}_{R,m}(t)-d_{R,m+1}(t))\hat{e}'P_{R,m}(t)\hat{e}, \nonumber \\
&\forall t\in J_{\nu},\nu\in\mathbb{N},\hat{e}\in\Rat{m} \label{Linmpl1:follows:from2}
\end{align}

\noindent \textbf{Establishment of \eqref{Linmpl1:follows:from2}, \eqref{map:PRm:properties} and \eqref{map:dRm:properties} for $m=m+1$:} We impose the following additional requirement for the candidate function $d_{R,m+1}(\cdot)$:
\begin{equation} \label{dRmpl1:lt:dRbm}
d_{R,m+1}(t)\le\bar{d}_{R,m}(t),\forall t\in J_{\nu},\nu\in\mathbb{N}.
\end{equation}

\noindent Then, by taking into account \eqref{dRmpl1:lt:dRbm}, it follows that, in order to show \eqref{Linmpl1:follows:from2}, it suffices to show that
\begin{align}
&p_{R,m+1}(t)a_{n-m}(t-\tau,y_{\tau}(t),u_{\tau}(t))+\phi_{R,m}(t)\le \bar{d}_{R,m}(t)-d_{R,m+1}(t) \nonumber \\
&\forall t\in J_{\nu},\nu\in\mathbb{N}, \label{Linmpl1:follows:from3}
\end{align}

\noindent for suitable functions $p_{R,m+1,1},p_{R,m}\in C^{1}(\cup_{\nu\in\mathbb{N}}J_{\nu};{\mathbb R})$ and $d_{R,m+1}\in C^{1}(\cup_{\nu\in\mathbb{N}}J_{\nu};{\mathbb R})$, in such a way that \eqref{map:PRm:properties}, \eqref{map:dRm:properties} hold with $m=m+1$, and in addition  $d_{R,m+1}(\cdot)$ satisfies \eqref{dRmpl1:lt:dRbm}. We proceed to the explicit construction of these functions.

\noindent \textbf{Construction of the mappings ${p_{R}(\cdot)}$ and ${d_{R,m+1}(\cdot)}$:} Let
\begin{subequations} \label{constants:Mmpl1:taumpl1}
	\begin{align}
	M_{m+1,\nu}:=&\max\left\{ |\bar{d}_{R,m}(t)|+\phi_{m}(t):t\in J_{\nu}\right\},\nu\in\mathbb{N} \label{constant:Mmpl1} \\
	\delta\varepsilon_{m,\nu}:=&\min\left\{\frac{1}{4nM_{m+1,\nu}},\frac{1}{4n^{2}},\frac{a_{\nu}-\varepsilon_{m,\nu}-\bar{a}_{\nu}}{2},\frac{\bar{b}_{\nu}-b_{\nu}-\varepsilon_{m,\nu}}{2}\right\} \nonumber \\
	\varepsilon_{m+1,\nu}:=&\varepsilon_{m,\nu}+\delta\varepsilon_{m+1,\nu} \label{constant:taukpl1}
	\end{align}
\end{subequations}

\noindent {and let $d_{R,m+1}\in C^1(\cup_{\nu\in\mathbb{N}}J_{\nu};\Rat{})$, $p_{R,m+1}\in C^1(\cup_{\nu\in\mathbb{N}}J_{\nu};\Rat{})$ in such a way that}
\begin{align}
d_{R,m+1}(t)&\left\lbrace\begin{array}{ll} \vspace{0.15cm}
:=\bar{d}_{R,m}(t)-\phi_{R,m}(t),  & t\in[\bar{a}_{\nu},a_{\nu}-\varepsilon_{m+1,\nu}],\nu\in\mathbb{N} \\ \vspace{0.15cm}
\in [\bar{d}_{R,m}(t)-\phi_{R,m}(t),\bar{d}_{R,m}(t)], & t\in [a_{\nu}-\varepsilon_{m+1,\nu},a_{\nu}-\varepsilon_{m,\nu}],\nu\in\mathbb{N} \\ \vspace{0.15cm}
:=\bar{d}_{R,m}(t), & t\in[a_{\nu}-\varepsilon_{m,\nu},b_{\nu}+\varepsilon_{m,\nu}],\nu\in\mathbb{N} \\ \vspace{0.15cm}
\in [\bar{d}_{R,m}(t)-\phi_{R,m}(t),\bar{d}_{R,m}(t)], & t\in [b_{\nu}+\varepsilon_{m,\nu},b_{\nu}+\varepsilon_{m+1,\nu}],\nu\in\mathbb{N} \\ \vspace{0.15cm}
:=\bar{d}_{R,m}(t)-\phi_{R,m}(t),  & t\in[b_{\nu}+\varepsilon_{m+1,\nu},\bar{b}_{\nu}],\nu\in\mathbb{N}
\end{array}\right. \label{map:dR2} \\
p_{R,m+1}(t)&:={\frac{\bar{d}_{R,m}(t)-\phi_{R,m}(t)-d_{R,m+1}(t)}{a_{n-m}(t-\tau,y_{\tau}(t),u_{\tau}(t))}},t\in J_{\nu},\nu\in\mathbb{N} \label{map:pR:keq2:dfn}
\end{align}

For $d_{R,m+1}(\cdot)$ and $p_{R,m+1}(\cdot)$ as defined in \eqref{map:dR2} and \eqref{map:pR:keq2:dfn} respectively, it can be seen that \eqref{Linmpl1:follows:from3} holds for all $t\in J_\nu$ and in addition, $d_{R,m+1}(t)\leq \bar{d}_{R,m}(t)$ for all $t\in J_{\nu}$. The latter, together with \eqref{map:dRbm:properties} and \eqref{map:dR2}, imply that
\begin{equation}\label{map:dRm+1:property1}
d_{R,m+1}(t)\leq \bar{d}_{R,m}(t)\leq -\bar{\sigma}_{R}(t)
\end{equation}

\noindent for all $t\in J_{\nu}\setminus[a_{\nu}-\varepsilon_{m+1,\nu},b_{\nu}+\varepsilon_{m+1,\nu}]$. {Notice that, due to the previous induction step and the result of Lemma \ref{lemma}, $\bar{d}_{R,m}(\cdot)$ and $\phi_{R,m}(\cdot)$ are strongly causal on each $J_{\nu}$, $\nu\in\mathbb{N}$, and therefore, $d_{R,m+1}(\cdot)$ and $p_{R,m+1}(\cdot)$ as defined in \eqref{map:dR2} and \eqref{map:pR:keq2:dfn} are also strongly causal on each $J_{\nu}$, $\nu\in\mathbb{N}$.}  We next show that  $d_{R,m+1}(\cdot)$ as defined in \eqref{map:dR2}, satisfies
\begin{equation*}
\int_{S}d_{R,m+1}(s)ds\geq -\frac{m+1}{n}
\end{equation*}

\noindent for all $S\subset [a_\nu-\varepsilon_{m+1,\nu},b_{\nu}+\varepsilon_{m+1,\nu}]\setminus I_\nu$. Indeed, from \eqref{map:dRbm:properties}, \eqref{constants:Mmpl1:taumpl1}, \eqref{map:dR2} and \eqref{map:dRm+1:property1} we get
\begin{align}
\int_{S}d_{R,m+1}(s)ds&=
\int_{([a_\nu-\varepsilon_{m+1,\nu},b_{\nu}+\varepsilon_{m+1,\nu}]\setminus [a_\nu-\varepsilon_{m,\nu},b_{\nu}+\varepsilon_{m,\nu}])\cap S}d_{R,m+1}(s)ds\nonumber\\
&+\int_{([a_\nu-\varepsilon_{m,\nu},b_{\nu}+\varepsilon_{m,\nu}]\setminus I_\nu)\cap S}d_{R,m+1}(s)ds\nonumber\\
\geq& -\int_{([a_\nu-\varepsilon_{m+1,\nu},b_{\nu}+\varepsilon_{m+1,\nu}]\setminus [a_\nu-\varepsilon_{m,\nu},b_{\nu}+\varepsilon_{m,\nu}])\cap S}(|\bar{d}_{R,m}(s)|+\phi(s))ds\nonumber\\
&+\int_{([a_\nu-\varepsilon_{m,\nu},b_{\nu}+\varepsilon_{m,\nu}]\setminus I_\nu)\cap S}\bar{d}_{R,m}(s)ds\nonumber\\
\geq&-2M_{m+1,\nu}\delta\varepsilon_{m,\nu}-\frac{m}{n}-\frac{1}{2n}\nonumber=-\frac{1}{2n}-\frac{2m+1}{2n}\nonumber\\
=&-\frac{m+1}{n},\label{map:dRm+1:prop2}
\end{align}

\noindent  We conclude, by taking into account \eqref{map:dRm+1:property1} and \eqref{map:dRm+1:prop2}, that all properties of \eqref{map:dRm:properties} are satisfied.

\noindent \textbf{Construction of the map ${p_{R,m+1,1}(\cdot)}$:}

It remains to determine a function $p_{R,m+1,1}\in C^1([t_0,\infty), \mathbb{R})$ such that \eqref{map:PRm:properties} holds for the general case, namely $P_{R,m+1}(t)>I_{(m+1)\times (m+1)}$. Define
\begin{equation}\label{map:pRm+1}
p_{R,m+1,1}(t):=\frac{p^2_{R,m+1}(t)\det(P_{R,m-1}(t)-I_{(m-1)\times (m-1)})}{\det(P_{R,m}(t)-I_{m\times m})}+L, \, \, t\geq t_0.
\end{equation}

\noindent Since $L>1$ and $P_{R,m}(t)>I_{m\times m}$, we deduce by using \eqref{map:pRm+1} that
\begin{align}
\det&(P_{R,m+1}(t)-I_{(m+1)\times (m+1)})=\nonumber\\
&=\det\left(\begin{matrix} p_{R,m+1,1}(t)-1 & \begin{matrix} p_{R,m+1}(t) & 0 & \cdots & 0 \end{matrix} \\ \begin{matrix} p_{R,m+1}(t) \\ 0 \\ \vdots \\ 0 \end{matrix} & \boxed{P_{R,m}(t)-I_{m\times m}} \end{matrix}\right)\nonumber\\
&=(L-1)\det(P_{R,m}(t)-I_{m\times m})>0, \, \, \forall t\geq t_0. \label{map:pRm+1:det}
\end{align}

\noindent Notice that for all $t\in J_\nu\setminus[a_{\nu}-\varepsilon_{m+1,\nu},b_\nu+\varepsilon_{m+1,\nu}]$ the mappings $p_{R,m+1}(\cdot)$ and $p_{R,m+1,1}(\cdot)$ as defined in \eqref{map:pR:keq2:dfn} and \eqref{map:pRm+1} respectively, satisfy $p_{R,m+1,1}(t)=L$ and $p_{R,m+1}(t)=0$. Therefore, all properties of \eqref{map:PRm:properties} hold.

It follows from \eqref{Lyapunov:inequality:aux:m}, \eqref{map:dRbm:properties}, \eqref{Linmpl1:follows:from2}, \eqref{Linmpl1:follows:from3}-\eqref{map:pR:keq2:dfn}, \eqref{map:dRm+1:prop2}, \eqref{map:pRm+1} and \eqref{map:pRm+1:det} that the induction hypothesis holds for $m:=m+1$. Therefore the proof of Claim 1 is complete.
\end{IEEEproof}
	
\noindent In order to show that \eqref{map:PR:properties}-\eqref{Lyapunov:ker:inequality} hold for $P_R(\cdot)$ and $d_R(\cdot)$ as defined in \eqref{map:PR:properties} and \eqref{map:dR:properties} respectively, we exploit Claim 1 with $m=n$. In particular, we consider the intervals $\{\mathcal{A}_{\nu}\}_{\nu\in\mathbb{N}}:=\{[\alpha_{\nu},\beta_{\nu}]\}_{\nu\in\mathbb{N}}$, given by {\eqref{epsilonm:intervals} with $m=n$, namely,}
\begin{equation}\label{intervals:Anu}
	\mathcal{A}_{\nu}:=[a_{\nu}-\varepsilon_{n,\nu},b_{\nu}+\varepsilon_{n,\nu}],\nu\in\mathbb{N}
\end{equation}

\noindent and define the mappings $P_{R}:[\bar{t}_0,\infty)\to\Rat{n\times n}$ and $d_{R}:[\bar{t}_0,\infty)\to\Rat{}$ as
\begin{equation} \label{map:PR:dfn}
	P_{R}(t):=\left\lbrace\begin{array}{ll}
	{\rm diag}\{L,\ldots,L\}, & t\in[\bar{t}_0,\infty)\setminus \cup_{\nu\in\mathbb{N}}\mathcal{A}_{\nu} \\
	P_{R,n}(t), & t\in\mathcal{A}_{\nu},\nu\in\mathbb{N}
	\end{array}\right.
\end{equation}
	
\noindent and
\begin{equation} \label{map:dR:dfn}
	d_{R}(t):=\left\lbrace\begin{array}{ll}
	-\bar{\sigma}_{R}(t), & t\in[\bar{t}_0,\infty)\setminus \cup_{\nu\in\mathbb{N}}\mathcal{A}_{\nu} \\
	d_{R,n}(t), & t\in\mathcal{A}_{\nu},\nu\in\mathbb{N}
	\end{array}\right.
\end{equation}
	
\noindent with $L(>1)$ as given above and $\mathcal{A}_{\nu}$, $P_{R,n}(\cdot)$, $d_{R,n}(\cdot)$ and $\bar{\sigma}_{R}(\cdot)$ as in \eqref{intervals:Anu}, \eqref{map:PRm:properties}, \eqref{map:dRm:properties} and \eqref{map:sigmaRbar}, respectively. It then follows {from \eqref{ais:ne0:inInu}-\eqref{map:PRm:properties} and \eqref{map:PR:dfn}} that $P_{R}\in C^1([\bar{t}_0,\infty);\Rat{n\times n})$ and that $d_{R}:[\bar{t}_0,\infty)\to\Rat{}$ is piecewise continuous. {Moreover, by invoking \eqref{intervals:Inu:Jnu}-\eqref{epsilonm:intervals} and \eqref{intervals:Anu}-\eqref{map:dR:dfn}, it follows that for each $\nu\in\mathbb{N}$, the mappings above are causal and {strongly} causal on each $\mathcal{A}_{\nu}$ with respect to $O(t_0,M)$}. In addition, it can be shown that $P_R(\cdot)$ and $d_R(\cdot)$, as defined by \eqref{map:PR:dfn} and \eqref{map:dR:dfn} respectively, satisfy \eqref{map:PR:properties}-\eqref{Lyapunov:ker:inequality}, with $\kappa_R(\cdot)$ as given by \eqref{map:kappaR:dfn} and
\begin{align}
	A(t,q,y,u):= & A_n(t,q,y,u) \label{map:A:dfn} \\
	H:= & H_n \label{map:H:dfn}
\end{align}
	
\noindent where $A_{n}(\cdot,\cdot,\cdot,\cdot)$ and $H_n$ are defined in \eqref{map:Am} and \eqref{map:Hm}, respectively. {We proceed with the establishment of \eqref{map:PR:properties}-\eqref{Lyapunov:ker:inequality}.}
	
\noindent \textbf{Establishment of \eqref{map:PR:properties}:} By taking into account the third property of \eqref{map:PRm:properties}, \eqref{map:PR:dfn} and the fact that due to \eqref{epsilonm:intervals} the left endpoint of $\mathcal{A}_1$ is greater than $t_0$, we deduce that both requirements of \eqref{map:PR:properties} are fulfilled.
	
\noindent \textbf{Establishment of \eqref{map:dR:properties}:} By taking into account \eqref{map:dRm:properties}, \eqref{intervals:Anu} and \eqref{map:dR:dfn} it follows that $d_{R}(\cdot)$ satisfies
\begin{equation} \label{map:dR:properties:aux}
	\begin{array}{ll} \vspace{0.15cm}
	d_R(t)=\frac{\int_{\bar{t}_{0}}^{T_{2}+\tau}\bar{\sigma}_{R}(s)ds+2}{b_{1}-a_{1}}, & \forall t\in I_{1} \\ \vspace{0.1cm}
	d_R(t)=\frac{\int_{T_{\nu}+\tau}^{T_{\nu+1}+\tau}\bar{\sigma}_{R}(s)ds+2}{b_{\nu}-a_{\nu}}, & \forall t\in I_{\nu},\nu=2,3,\ldots \\ \vspace{0.15cm}
	\int_{S}d_R(s)ds\ge -1, & \forall S\subset [a_{\nu}-\varepsilon_{n,\nu},b_{\nu}+\varepsilon_{n,\nu}]\setminus I_{\nu},\nu=1,2,\ldots \\
	\end{array}
\end{equation}
	
\noindent In order to prove the desired \eqref{map:dR:properties}, we consider two cases.
	
\noindent \textbf{Case A: $t\in[\bar{t}_0,T_1+\tau)$.} Then, it follows from \eqref{map:dR:dfn}, \eqref{map:dR:properties:aux} and \eqref{map:kappaR:dfn} that
\begin{align*}
	\int_{[\bar{t}_0,t)}d_R(s)ds & =\int_{[\bar{t}_0,t)\cap [\bar{t}_0,T_1+\tau]}d_R(s)ds=\int_{[\bar{t}_0,t)\cap ([\bar{t}_0,T_1+\tau]\setminus\mathcal{A}_1)}d_R(s)ds  \\
	& + \int_{[\bar{t}_0,t)\cap (\mathcal{A}_1\setminus I_1)}d_R(s)ds+\int_{[\bar{t}_0,t)\cap I_1}d_R(s)ds\ge\int_{[\bar{t}_0,T_1+\tau]\setminus\mathcal{A}_1}d_R(s)ds-1 \\
	& > -\int_{[\bar{t}_0,T_1+\tau]}\bar{\sigma}_R(s-\tau)ds-1\ge\kappa_R(t)
\end{align*}
	
\noindent \textbf{Case B: $t\in[T_{\nu}+\tau,T_{\nu+1}+\tau), \nu\in\mathbb{N}$.} Then we get that
\begin{equation} \label{integral:decomposition}
	\int_{[\bar{t}_0,t)}d_R(s)ds=\sum_{i=1}^{\nu}\int_{[T_{\nu-1}+\tau,T_{\nu}+\tau)}d_R(s)ds+\int_{[T_{\nu}+\tau,t]}d_R(s)ds
\end{equation}
	
\noindent By performing similar manipulations with those in Case A, we obtain that the second term in \eqref{integral:decomposition} satisfies
\begin{equation} \label{integral:decomposition:term2}
	\int_{[T_{\nu}+\tau,t]}d_R(s)ds>-\int_{[T_{\nu}+\tau,T_{\nu+1}+\tau)}\bar{\sigma}_R(s-\tau)ds-1
\end{equation}
	
\noindent For the first term, by exploiting \eqref{map:dR:dfn} and \eqref{map:dR:properties:aux}  we deduce that
\begin{align}
	\sum_{i=1}^{\nu} & \int_{[T_{i-1}+\tau,T_{i}+\tau)}d_R(s)ds =\sum_{i=1}^{\nu}\int_{[T_{i-1}+\tau,T_i+\tau)\setminus\mathcal{A}_i}d_R(s)ds+\sum_{i=1}^{\nu} \int_{\mathcal{A}_i\setminus I_i}d_R(s)ds \nonumber \\
	& + \sum_{i=1}^{\nu}\int_{I_i}d_R(s)ds\ge\sum_{i=1}^{\nu}\left(-\int_{[T_{i-1}+\tau,T_i+\tau)\setminus\mathcal{A}_i}\bar{\sigma}_R(s-\tau)ds\right)+\sum_{i=1}^{\nu}(-1) \nonumber \\
	& + (b_1-a_1)\frac{\int_{\bar{t}_{0}}^{T_{2}+\tau}\bar{\sigma}_{R}(s)ds+2}{b_{1}-a_{1}}+(b_{\nu}-a_{\nu})\sum_{i=2}^{\nu}\frac{\int_{T_{\nu}+\tau}^{T_{\nu+1}+\tau}\bar{\sigma}_{R}(s)ds+2}{b_{\nu}-a_{\nu}} \nonumber \\
	 \ge& -\int_{[\bar{t}_0,T_{\nu}+\tau)}\bar{\sigma}_R(s-\tau)ds-\nu+2\nu+\int_{[\bar{t}_0,T_{\nu+1}+\tau)}\bar{\sigma}_R(s-\tau)ds \nonumber \\
	 =&\,\, \nu+\int_{[T_{\nu}+\tau,T_{\nu+1}+\tau)}\bar{\sigma}_R(s-\tau)ds. \label{integral:decomposition:term1}
\end{align}
	
\noindent Hence, we get from \eqref{map:kappaR:dfn}, \eqref{integral:decomposition}, \eqref{integral:decomposition:term2}  and \eqref{integral:decomposition:term1} that $\int_{[\bar{t}_0,t)}d_R(s)ds>\nu-1\ge\kappa_R(t)$ and conclude that \eqref{map:dR:properties} holds in this case as well.
	
\noindent \textbf{Establishment of \eqref{Lyapunov:ker:inequality}:} We consider again two cases.
	
\noindent \textbf{Case A: $t\in\mathcal{A}_{\nu}$ for some $\nu\in\mathbb{N}$.} Then, it follows from \eqref{map:PRm:properties}, \eqref{map:dRm:properties}, \eqref{Lyapunov:inequality:m}, \eqref{map:PR:dfn} and  \eqref{map:dR:dfn} that \eqref{Lyapunov:ker:inequality} is satisfied with $A(\cdot,\cdot,\cdot,\cdot)$ and $H$ as given by \eqref{map:A:dfn}, \eqref{map:H:dfn}.

\noindent \textbf{Case B: $ t\in[\bar{t}_0,\infty)\setminus \cup_{\nu\in\mathbb{N}}\mathcal{A}_{\nu}$.} In this case, we obtain from \eqref{map:A:triang:Frobenius:norm}, \eqref{map:PR:dfn} and \eqref{map:dR:dfn} that for all $e\in\ker H$ and $q\in Q_{R}(t-\tau)$ it holds
\begin{align*}
	e'P_{R}(t)A(t-\tau,q,y_{\tau}(t))e+\frac{1}{2}e'\dot{P}_{R}(t)e & =e'P_{R}(t)A(t-\tau,q,y_{\tau}(t),u_{\tau}(t))e \\
	& \le L\bar{\sigma}_{R}(t)|e|^2=-d_R(t)e'P_{R}(t)e.
\end{align*}
Thus, Claim 1 guarantees that \eqref{map:PR:properties}-\eqref{Lyapunov:ker:inequality} of A2 are fulfilled.
	
We conclude that for every $R>0$ with $B_{R}\cap M\neq\emptyset$, both Hypothesis \ref{hypothesis:det} and \ref{hypothesis:switchin} hold, thus, according to Proposition \ref{Proposition:switching} the IDSODP is solvable for \eqref{system:triangular} with respect to $(M,\U)$. Statement (ii) of Proposition \ref{Theorem} is a direct consequence of Claim 1, Proposition \ref{Proposition:state_deter} and the fact that the initial states of \eqref{system:triangular} belong to the intersection of $M$ with a given ball $B_{R}$ of radius $R>0$.
\end{IEEEproof}

\section{Conclusions}
In this paper, sufficient conditions are established for the solvability of the observer design problem for a class of nonlinear triangular control systems. The Luenberger-type observer we propose is in general time-varying and the state estimation is achieved with an arbitrarily small delay. The global state determination is treated by a switching observer methodology. The main result is based on a forwarding inductive procedure which extends the approach employed in \cite{BdTj13a}.

\renewcommand{\theequation}{A.\arabic{equation}}
\appendix
\begin{IEEEproof}[Proof of Lemma \ref{lemma} (Outline)]
Let $t_{0}$, $\tau$, $a$, $b$, $y(\cdot)$, $u(\cdot)$ and consider a triple of mappings $d(\cdot)$, $P(\cdot)$ and $\bar{d}(\cdot)$ as given in the statement of the lemma. We proceed with the construction of the function $\phi(\cdot)$ on $[a,b]$ and define for each $t\in[a,b]$, $q\in\Rat{\ell}$ and $e\in\Rat{r}$ the mappings
\begin{align}
D(t,q,e):=&e'P(t)A(t-\tau,q,y_{\tau}(t),u_{\tau}(t))e+\frac{1}{2}e'\dot{P}(t)e+\bar{d}(t)e'P(t)e \label{map:D} \\
K(t):=&\{w\in\Rat{r}:|w|=1,D(t,q,w)<0,\forall q\in Q(t-\tau)\}. \label{set:K}
\end{align}

\noindent Also, continuity of $y(\cdot)$, $u(\cdot)$ and the mappings involved in the right hand side of \eqref{map:D}, imply that $D(\cdot,\cdot,\cdot)$ is continuous. Notice that due to \eqref{rank:H}, \eqref{dbar:lt:d} and \eqref{Lyapunov:inequality:aux} the set $K(t)$ is nonempty, since it includes all vectors $w\in\Rat{r}$ with $|w|=1$ and $w\in{\rm ker}H(t-\tau,u_{\tau}(t))\ne\emptyset$. Indeed, let $w\in\Rat{r}$ with $|w|=1$ and $w\in{\rm ker}H(t-\tau,u_{\tau}(t))$. Then, by using  \eqref{dbar:lt:d}, \eqref{Lyapunov:inequality:aux} and by taking into account that $P(\cdot)$ is positive definite, we deduce that $D(t,q,w)\le(\bar{d}(t)-d(t))w'P(t)w<0$ for all $q\in Q(t-\tau)$ and hence that $w\in K(t)$ which asserts that $K(t)\ne\emptyset$. Thus we established the implication
\begin{equation} \label{set:K:property}
w\in{\rm ker} H(t-\tau,u_{\tau}(t))\;{\rm and}\; |w|=1\Rightarrow w\in K(t).
\end{equation}

\noindent In the sequel, for each $t\in [a,b]$, we adopt the notation $K^{c}(t)$ to indicate the complement of $K(t)$ with respect to the unit sphere in $\Rat{n}$, namely, $K^{c}(t):=\{w\in\Rat{r}:|w|=1,w\notin K(t)\}$. Hence, we get from \eqref{set:K:property} that
\begin{equation} \label{set:Kc}
K^{c}(t)=\{w\in\Rat{r}:|w|=1\;{\rm and}\;D(t,q,w)\ge 0,\;\textup{for some}\; q\in Q(t-\tau)\},
\end{equation}

\noindent {and similarly to the proof of \cite[Lemma 2.1]{BdTj13a}, it can be shown, by exploiting \eqref{set:Kc}, the CP property and continuity of $D(t,\cdot,\cdot)$, that for every $t\in[a,b]$ the set $K^{c}(t)$ is closed.} Next, we consider the map $\omega:[a,b]\to [0,\infty]$ defined as
\begin{equation}\label{map:omega}
\omega(t):=\left\lbrace \begin{array}{ll} \vspace{0.15cm}
\min\{|H(t-\tau,u_{\tau}(t))w|:w\in K^{c}(t)\}, & {\rm if}\; K^{c}(t)\ne\emptyset \\ \vspace{0.15cm}
\infty, & {\rm if}\; K^{c}(t)=\emptyset
\end{array}\right.
\end{equation}

\noindent Notice that for every $t\in [a,b]$ the set $\{|H(t-\tau,u_{\tau}(t))w|:w\in K^{c}(t)\}$ is compact, whenever $K^{c}(t)\ne\emptyset$ and hence $\omega(\cdot)$ is well defined and satisfies $\omega(t)>0$ for all $t\in [a,b]$. It also holds
\begin{equation} \label{map:omega:property}
\inf\{\omega(t):t\in [a,b]\}>0.
\end{equation}

\noindent {The proof of \eqref{map:omega:property} is quite similar to that given in \cite{BdTj13a} and is omitted.} Next, we define the mapping $\bar{\omega}:[a,b]\to \RgeO$ as
\begin{equation}\label{map:omega:bar}
\bar{\omega}(t):=\left\lbrace \begin{array}{ll} \vspace{0.15cm}
\frac{1}{\omega^{2}(t)} & {\rm if}\; K^{c}(t)\ne\emptyset \\ \vspace{0.15cm}
0 & {\rm if}\; K^{c}(t)=\emptyset
\end{array}\right.
\end{equation}

\noindent It then follows from \eqref{map:omega:property} and \eqref{map:omega:bar} that there exists a constant $M>0$ such that $\sup\{\bar{\omega}(t):t\in[a,b]\}\le M.$ Also, define for $t\in[a,b]$
\begin{equation} \label{map:C}
C(t):=\sup\left\lbrace\bar{\omega}(t)\left(|P(t)||A(t-\tau,q,y_{\tau}(t),u_{\tau}(t))|+\frac{1}{2}|\dot{P}(t)|+|\bar{d}(t)||P(t)|\right):q\in Q(t-\tau)\right\rbrace.
\end{equation}

\noindent  It is then straightforward to construct a function $\phi\in C^{1}([a,b];{\mathbb R}_{\ge 0})$ satisfying
\begin{equation} \label{map:phi}
\phi(t)>C(t), \forall t\in[a,b].
\end{equation}

\noindent  From the above constructions, the hypothesis that $P(\cdot)$, $d(\cdot)$ and $\bar{d}(\cdot)$ are {strongly}  causal on $[a,b]$ with respect to $\Omega(t_{0};W)$ and our assumption that $b-a<\tau$, it follows that $\phi(\cdot)$ is also {strongly} causal on $[a,b]$ with respect to $\Omega(t_{0};W)$. Next, notice that the desired \eqref{Lyapunov:inequality:aux} is equivalent to
\begin{align}
&w'P(t)A(t-\tau,q,y_{\tau}(t),u_{\tau}(t))w+\tfrac{1}{2}w'\dot{P}(t)w\le\phi(t)|H(t-\tau,u_{\tau}(t))w|^{2}-\bar{d}(t)w'P(t)w, \nonumber \\
& \forall t\in [a,b],w\in\Rat{r}:|w|=1,q\in Q(t-\tau) \label{Lyapunov:inequality:aux2}
\end{align}

\noindent hence, in order to prove \eqref{Lyapunov:inequality:aux}  it suffices to show that \eqref{Lyapunov:inequality:aux2} is fulfilled. Indeed, let $t\in [a,b]$ and $w\in K(t)$. Then the desired \eqref{Lyapunov:inequality:aux2} is a consequence of \eqref{map:D}, \eqref{set:K} and the fact that $\phi(t)> 0$. Finally, if $K^{c}(t)\ne\emptyset$ and $w\in K^{c}(t)$, then in order to show \eqref{Lyapunov:inequality:aux2}, it suffices to show that $
\sup\left\lbrace |P(t)||A(t-\tau,q,y_{\tau}(t),u_{\tau}(t))|+\frac{1}{2}|\dot{P}(t)|+|d(t)||P(t)|:q\in Q(t-\tau)\right\rbrace\le\phi(t)\omega^{2}(t)$ which is a consequence of \eqref{map:omega:bar}, \eqref{map:C} and \eqref{map:phi}. This completes the proof of Lemma \ref{lemma}.
\end{IEEEproof}

\end{document}